\documentclass[10pt]{article}
\usepackage{amsfonts}
\usepackage[leqno]{amsmath}
\usepackage{graphicx}
\usepackage{latexsym}
\usepackage{amsmath,amsfonts,amssymb,amsthm,mathrsfs,euscript,makeidx,color}
\usepackage{enumerate}
\allowdisplaybreaks[1]
\usepackage{multirow}

\def\sqr#1#2{{\vcenter{\vbox{\hrule height.#2pt
              \hbox{\vrule width.#2pt height#1pt \kern#1pt \vrule width.#2pt}
              \hrule height.#2pt}}}}
\def\signed #1{{\unskip\nobreak\hfil\penalty50
              \hskip2em\hbox{}\nobreak\hfil#1
              \parfillskip=0pt \finalhyphendemerits=0 \par}}
\def\endpf{\signed {$\sqr69$}}

\def\5n{\negthinspace \negthinspace \negthinspace \negthinspace \negthinspace }
\def\4n{\negthinspace \negthinspace \negthinspace \negthinspace }
\def\3n{\negthinspace \negthinspace \negthinspace }
\def\2n{\negthinspace \negthinspace }
\def\1n{\negthinspace }

\def\dbC{\mathbb{C}}

\def\dbH{\mathbb{H}}

\def\dbR{\mathbb{R}}
\def\dbS{\mathbb{S}}

\def\sR{\mathscr{R}}

\def\sU{\mathscr{U}}

\def\cO{{\cal O}}

\def\BN{{\bf N}}

\def\Be{{\bf e}}

\def\Bv{{\bf v}}

\def\Re{{\mathop{\rm Re}\,}}

\def\ds{\displaystyle}

\def\ns{\noalign{\ss}}
\def\no{\noindent}

\def\ss{\smallskip}
\def\ms{\medskip}

\def\q{\quad}
\def\qq{\qquad}
\def\hb{\hbox}

\def\({\Big (}
\def\){\Big )}
\def\[{\Big[}
\def\]{\Big]}

\def\lan{\langle}
\def\ran{\rangle}

\def\rf{\eqref}

\def\a{\alpha}
\def\b{\beta}

\def\d{\delta}
\def\e{\varepsilon}

\def\k{\kappa}
\def\l{\lambda}
\def\m{\mu}
\def\n{\nu}
\def\si{\sigma}
\def\t{\tau}
\def\f{\varphi}
\def\th{\theta}

\def\i{\infty}

\def\Th{\Theta}

\def\F{\Phi}

\def\limsup{\mathop{\overline{\rm lim}}}
\def\liminf{\mathop{\underline{\rm lim}}}

\def\Ra{\mathop{\Rightarrow}}

\def\span{\mathop{\rm span}}

\def\h{\widehat}
\def\wt{\widetilde}

\def\cd{\cdot}
\def\cds{\cdots}

\def\deq{\triangleq}
\def\les{\leqslant}
\def\ges{\geqslant}
\def\pa{\partial}
\def\rf{\eqref}
\def\bde{\begin{definition}\label}
\def\ede{\end{definition}}
\def\be{\begin{equation}}
\def\bel{\begin{equation}\label}
\def\ee{\end{equation}}
\def\bt{\begin{theorem}\label}
\def\et{\end{theorem}}
\def\bc{\begin{corollary}\label}
\def\ec{\end{corollary}}
\def\bl{\begin{lemma}\label}
\def\el{\end{lemma}}
\def\bp{\begin{proposition}\label}
\def\ep{\end{proposition}}
\def\bas{\begin{assumption}\label}
\def\eas{\end{assumption}}
\def\br{\begin{remark}\label}
\def\er{\end{remark}}
\def\bex{\begin{example}\label}
\def\ex{\end{example}}
\def\ba{\begin{array}}
\def\ea{\end{array}}
\def\ben{\begin{enumerate}}
\def\een{\end{enumerate}}

\def\square#1{\vbox{\hrule\hbox{\vrule height#1%
     \kern#1\vrule}\hrule}}
\def\rectangle#1#2{\vbox{\hrule\hbox{\vrule height#1%
     \kern#2\vrule}\hrule}}

\font\tenbb=msbm10 \font\sevenbb=msbm7 \font\fivebb=msbm5

\newfam\bbfam
\scriptscriptfont\bbfam=\fivebb \textfont\bbfam=\tenbb
\scriptfont\bbfam=\sevenbb

\newtheorem{theorem}{\indent Theorem}[section]
\newtheorem{definition}[theorem]{\indent Definition}
\newtheorem{proposition}[theorem]{\indent Proposition}
\newtheorem{corollary}[theorem]{\indent Corollary}
\newtheorem{lemma}[theorem]{\indent Lemma}
\newtheorem{remark}[theorem]{\indent Remark}
\newtheorem{example}[theorem]{\indent Example}

\newtheorem{assumption}[theorem]{\indent Assumption}

\makeatletter
   
   \@addtoreset{equation}{section}
\makeatother

\begin{document}
\title{\bf Infinite Horizon Linear Quadratic Overtaking\\ Optimal Control Problems}

\author{Jianping Huang\footnote{School of Mathematics and Statistics, Central South University, Changsha, Hunan, 410083, P.R. China (huangjianping@csu.edu.cn). This author was supported in part by Hunan Provincial Innovation Foundation for Postgraduate Grant CX20190093.},~~~Jiongmin Yong\footnote{Department of
Mathematics, University of Central Florida, Orlando, FL 32816, USA (Jiongmin.Yong@ucf.edu). This author was supported in part by NSF Grant DMS-1812921.},~~~Huacheng Zhou\footnote{School of Mathematics and Statistics, Central South University, Changsha, Hunan, 410083, P.R. China (hczhou@amss.ac.cn). This author was supported in part by the National Natural Science Foundation of China  Grant 61803386.}}

\maketitle

\begin{abstract}
A linear control system with quadratic cost functional over infinite time horizon is considered without assuming controllability/stabilizability condition and the global integrability condition for the nonhomogeneous term of the state equation and the weight functions in the linear terms in the running cost rate function. Classical approaches do not apply for such kind of problems. Existence and non-existence of overtaking optimal controls in various cases are established. Some concrete examples are presented. These results show that the overtaking optimality approach can be used to solve some of the above-mentioned problems and at the same time, the limitation of this approach is also revealed.

\end{abstract}

\bf Keywords. \rm linear quadratic problem, overtaking optimal control, controllability.

\ms

\bf AMS 2020 Mathematics Subject Classification. \rm 49J15, 49N10, 93B05

\section{Introduction}

Investigation of infinite time horizon optimal control problems can be traced back to the work of Ramsey in 1928 on {\sl a mathematical theory of saving} \cite{Ramsey 1928}. There is a big number of follow-up works, for examples, von Weizs\"acker \cite{Weizsacker 1965}, Arrow \cite{Arrow 1968}, Arrow--Kurtz \cite{Arrow-Kurtz 1970}, Halkin \cite{Halkin 1974}, Brock--Haurie \cite{Brock-Haurie 1976}, to mention a few for the period of 1960--1970s, and there were many more afterwards. For general (nonlinear) continuous-time controlled dynamics with the performance (cost/payoff) functional in infinite time horizons, to treat the situation that the performance functional is possibly not well-defined over the infinite time horizon, von Weizs\"acker introduced the so-called {\it overtaking optimization} approach in 1965  (\cite{Weizsacker 1965}), which, ``approximately'' compares the values of the performance functional over every finite interval. See \cite{Samuelson-Moussavian 1985, Carlson-Haurie-Leizarowitz 1991, Tan-Rugh 1998}, and references cited therein. We will make this precise later in the current paper. There are some other relevant works on this class of problems, without using overtaking optimality, see \cite{Aseev-Kryazhimskii 2007, Basco-Cannarsa-Frankowska 2018, Buckdahn-Li-Quincampoix-Renault 2020} and rich references cited therein. On the other hand, standard linear-quadratic optimal control problem (LQ problem, for short) in infinite time horizons is well-understood (\cite{Belyakov 2019, Kalman 1960, Willems 1971, Wonham 1985, Anderson-Moore 1989, Huang-Li-Yong 2015, Skritek-Veliov-2015, Sun-Yong 2018, Sun-Yong 2020}). However, we still find some interesting and challenging LQ problems relevant to the overtaking optimality. To elaborate that, let us begin with the following controlled linear ordinary differential equation:
\bel{state1}\left\{\2n\ba{ll}
\ds\dot X(s)=AX(s)+Bu(s)+b(s),\qq s\in[t,\i),\\
\ns\ds X(t)=x,\ea\right.\ee
where $A\in\dbR^{n\times n}$ and $B\in\dbR^{n\times m}$ are called the {\it coefficients}, $b:[0,\infty)\to\dbR^n$, a locally integrable over $[0,\i)$, is called the {\it nonhomogeneous term}. Here, $\dbR^{m\times n}$ is the set of all $(m\times n)$ matrices, and $\dbR^n=\dbR^{n\times1}$. Then
for any {\it initial pair} $(t,x)\in[0,\i)\times\dbR^n$, and any {\it control} $u(\cd)\in\sU_{loc}[t,\i)$ with
\bel{sU}\sU_{loc}[t,\i)\equiv L^2_{loc}(t,\i;\dbR^m)\deq\Big\{u:[t,\i)\to\dbR^m\bigm|\int_t^T|u(s)|^2ds<\i,\q\forall T>t\Big\}.\ee
{\it state equation} \rf{state1} admits a unique solution $X(\cd)\equiv X(\cd\,;t,x,u(\cd))$ which is called the {\it state trajectory}.
To measure the performance of the control $u(\cd)$, we introduce the following {\it running cost rate function}
\bel{cost rate}g(s,x,u)\1n=\1n\lan Qx,x\ran\1n+\1n2\lan Sx,u\ran\1n+\1n\lan Ru,u\ran\1n+\1n2\lan q(s),x\ran\1n+\1n2\lan\rho(s),u\ran,\q(s,x,u)\1n\in\1n[0,\i)\1n\times\1n\dbR^n\1n
\times\1n\dbR^m,\ee
with $Q\in\dbS^n$, $S\in\dbR^{m\times n}$ and $R\in\dbS^m$ being some constant
matrices (called {\it quadratic weighting matrices}), and $q:[0,\infty)\to\dbR^n$, $\rho:[0,\infty)\to\dbR^m$ being some locally integrable functions (called {\it linear weighting functions}). Here, $\dbS^n$ is the set of all $(n\times n)$ symmetric matrices. Unlike the classical situation, we do not assume the {\it stabilizability} of system $[A,B]$ and functions $b(\cd)$, $q(\cd)$ and $\rho(\cd)$ are only assumed to be (square) integrable on each finite interval $[0,T]$. Formally, the running cost over any infinite time interval $[t,\i)$ reads
\bel{cost[0,i]}J(t,x;u(\cd))\equiv J_\i(t,x;u(\cd))=\int_t^\i g(s,X(s),u(s))ds.\ee
Clearly, for any $(t,x,u(\cd))\in[0,\i)\times\dbR^n\times\sU_{loc}[0,\i)$,
$J(t,x;u(\cd))$ might not be well-defined. Therefore, we define
\bel{U_ad^x}\sU^x_{\1n J}[t,\i)=\Big\{u(\cd)\in\sU_{loc}[t,\i)\bigm|J(t,x;u(\cd))\hb{ is well-defined}\Big\}.\ee
Then one can formulate the following LQ problem on $[0,\i)$.

\ms

{\bf Problem (LQ)$_\i$.} For any $(t,x)\in[0,\i)\times\dbR^n$, find a $\bar u(\cd)\in\sU^x_{\1n J}[t,\i)$ such that
\bel{J(bar u)*}J(t,x;\bar u(\cd))=\inf_{u(\cd)\in\sU^x_{\1n J}[t,\i)}J(t,x;u(\cd))
=V_\i(t,x).\ee

If $\bar u(\cd)\in\sU^x_{\1n J}[t,\i)$ satisfies \rf{J(bar u)*}, we call it an {\it open-loop optimal control}, the corresponding $\bar X(\cd)\equiv X(\cd\,;t,x,\bar u(\cd))$ and
$(\bar X(\cd),\bar u(\cd))$ are called an {\it open-loop optimal trajectory}, and an
{\it open-loop optimal pair}, respectively, for Problem (LQ)$_\i$. Also,
$V_\i(\cd\,,\cd)$ is called the value function of the problem (\cite{Sun-Yong 2020}).

\ms

Note that for Problem (LQ)$_\i$, requiring the cost functional $J(t,x;u(\cd))$ to be finite, it roughly implies that the running cost rate $g(s,X(s),u(s))$ approaches to zero as $s\to\infty$. In some applications, this might not be expected. For example, if there exists a persistent part of running cost, by which we mean that $g(s,X(s),u(s))$ has a positive lower bound. Such a situation happens if we consider the cost of some production process, as time goes by, due to the demand-driven production level and possible increase of the prices of raw material, cost of manpower, etc., one could not expect to have a decreasing cost rate. Another possible situation is some kind of approximate (not necessarily exact) seasonal impact leading to the fluctuation/near periodicity of the running cost rate. Hence, we should allow the cost functional $J(t,x;u(\cd))$ to be divergent. Mathematically, we should not assume the stabilizability condition for system $[A,B]$. Likewise, the global integrability of $b(\cd)$, $q(\cd)$ and $\rho(\cd)$ should not be assumed, either. Consequently, although the running cost rate $g(\cd\,,X(\cd),u(\cd))$ is integrable over any finite interval $[t,T]$, it is not necessarily integrable over $[t,\infty)$. As a matter of fact, Ramsey's original problem (\cite{Ramsey 1928}) already has such a feature. In the current paper, we are going to  explore the LQ problem in infinite time horizons without assuming the stabilizability of the system $[A,B]$ and the global integrability of $b(\cd),q(\cd),\rho(\cd)$ by means of overtaking optimality. For convenience, we call the problem without assuming the stabilizability of $[A,B]$ and the global integrability of $b(\cd)$, $q(\cd)$, $\rho(\cd)$ Problem (LQ) (to distinguish from Problem (LQ)$_\i$).

\ms

Now, let us briefly highlight our main ideas and list the main contributions of the current paper.

\ms

Let $\dbH_0=\span\big\{\sR(A^kB)\bigm|k\ges0\big\}$ which is the controllable subspace of system $[A,B]$. Let $\Pi:\dbR^n\to\dbH_0$ be the orthogonal projection and $\Pi^\perp=I-\Pi$. Decompose the state $X(\cd)$ as follows:
\bel{X=X+X}X(\cd)=\Pi X(\cd)+\Pi^\perp X(\cd)\equiv X_\Pi(\cd)+X_{\Pi^\perp}(\cd).\ee
Then the control only affects $X_\Pi(\cd)$, and does not affect $X_{\Pi^\perp}(\cd)$. Having such a decomposition, we can summarize our main results:

\ms

(i) If the equation for $X_\Pi(\cd)$ does not contain $X_{\Pi^\perp}(\cd)$, and in the running cost rate function, $(X_\Pi(\cd),u(\cd))$ and $X_{\Pi^\perp}(\cd)$ are separated; in addition, the LQ problem involving $(X_\Pi(\cd),u(\cd))$ admits an optimal control,
then Problem (LQ) admits an overtaking optimal control.

\ms

(ii) If the conditions of (i) fail, then under some mild conditions, the general Problem (LQ) can be reduced to the following special case: the system $[A,B]$ is controllable with $A$ stable and $b(\cd)=0$; in the running cost rate function, $\rho(\cd)=0$ and $q(\cd)$ might be just locally integrable. For such a special case, under various conditions, we will show Problem (LQ) does not admit overtaking optimal control if the control set is a linear subspace of $L^2(t,\i;\dbR^m)$.

\ms

(iii) In the case that overtaking optimal control does not exist over the linear space of controls, we will discuss the problem over some convex set of controls. Under certain conditions, one can obtain the existence of overtaking optimal controls.

\ms

The above results give a general picture of overtaking optimality approach for the infinite horizon LQ problems. They roughly show both the power and the limitation of the approach.

\ms

The rest of the paper is organized as follows. In Section 2, we present some preliminary results, including mathematical motivation of introducing overtaking optimality, and some standard results of LQ problems in infinite time horizons. Section 3 is devoted to a presentation of the best case for which the problem has overtaking optimal controls. In Section 4, we present some general results on overtaking optimal controls for Problem (LQ), and introduce a further reduction. Then in Section 5, we carefully discuss the non-existence of overtaking optimal controls for various situations. This actually shows the limitation of overtaking optimality approach (at least) to the LQ problems. A set of sufficient conditions for the existence of overtaking optimal controls is presented in Section 6. Finally, some concluding remarks are collected in Section 7.

\section{Preliminary Results}

For $t\in[0,\i)$, $p\ges1$, and Euclidean space $\dbH$ (say, $\dbR^n$, $\dbR^{m\times n}$, etc.), we define
$$\ba{ll}
\ns\ds L^p_{loc}(t,\infty;\dbH)=\Big\{\f:[t,\infty)\to\dbH\bigm|\int_t^T|\f(t)|^pdt<\infty,\q\forall\ T>t\Big\},\\
\ns\ds L^p(t,\infty;\dbH)=\Big\{\f:[t,\infty)\to\dbH\bigm|\int_t^\infty|\f(t)|^pdt<\infty\Big\},\\
\ns\ds C([t,\infty);\dbH)=\Big\{\f:[t,\infty)\to\dbH\bigm|\f(\cd)\hb{ is continuous}\Big\}.\ea$$
According to the above, we have $\sU_{loc}[t,\infty)=L^2_{loc}(t,\i;\dbR^m)$. We denote (comparing with $\sU_{\1n J}^x[t,\i)$ defined in \rf{U_ad^x})
\bel{cU_0}\sU_0[t,\infty)=L^2(t,\i;\dbR^m),\qq\sU_{ad}^x[t,\infty)=\Big\{u(\cd)\in\sU_0[t,\infty)\bigm|X(\cd\,;x,u(\cd))\in L^2(t,\infty;\dbR^n)\Big\}.\ee
We may likewise define $L^p(t,T;\dbH)$ and $C([t,T];\dbH)$. The following proposition collects some standard results for Problem (LQ)$_\i$. See \cite{Sun-Yong 2018} for details.

\bp{LQ results} \sl Let $[A,B]$ be {\it stabilizable}, i.e., there exists a $\Th\in\dbR^{m\times n}$, called a {\it stabilizer} of the system, such that
\bel{stabilizable}\si(A+B\Th)\subseteq\dbC^-\equiv\{\a+i\b\bigm|\a,\b\in\dbR,~\a<0\},\ee
where $\si(A+B\Th)$ is the set of all eigenvalues of $A+B\Th$. Also, the function $b(\cd)\in L^1(0,\infty;\dbR^n)$. Let $Q\in\dbS^n$, $S\in\dbR^{m\times n}$, and $R\in\dbS^m$ satisfy
\bel{QSR}R>0,\q Q-S^\top R^{-1}S\ges0,\ee
and let $q(\cd)\in L^1(0,\infty;\dbR^n)$, $\rho(\cd)\in L^2(0,\infty;\dbR^m)$. Then for each $(t,x)\in[0,\i)\times\dbR^n$,
\bel{U non-empty}\varnothing\ne\sU^x_{ad}[t,\i)\subseteq\sU_{\1n J}^x[t,\i),\qq\forall(t,x)\in[0,\i)\times\dbR^n,\ee
and Problem {\rm(LQ)$_\i$} admits a uniquely open-loop control $\bar u(\cd)$. Moreover, the following algebraic Riccati equation
\bel{ARE}PA+A^\top P-(B^\top P+S)^\top R^{-1}(B^\top P+S)+Q=0,\ee
admits a solution $P>0$ such that $\si\big(A-B R^{-1}(B^\top P+S)\big)\subseteq\dbC^-$,
and the following equation:
$$\left\{\2n\ba{ll}
\ds\dot\eta(s)=-\big[A^\top-(PB+S^\top)R^{-1}B^\top\big]\eta(s)+(PB+S^\top)R^{-1}\rho(s)
-Pb(s)-q(s),\q s\in[t,\infty),\\
\ns\ds\eta(\infty)=0,\ea\right.$$
admits a unique solution $\eta(\cd)\in L^2(t,\i;\dbR^n)$. The unique open-loop optimal control $\bar u(\cd)$ admits the following {\it closed-loop} representation:
$$\bar u(s)=\bar\Th\bar X(s)+\bar v(s),\qq s\in[t,\infty),$$
where $\bar\Th=-R^{-1}(S+B^{\top}P)$ and $\bar v(\cd)=-R^{-1}[B^{\top}\eta(\cd)+\rho(\cd)]$.  Moreover, the value function is given by
$$V_\i(t,x)=\lan Px,x\ran+2\lan\eta(t),x\ran+\int_t^\infty\[2\lan\eta(s),b(s)\ran-\lan R^{-1}[B^\top\eta(s)+\rho(s)],B^\top\eta(s)+\rho(s)\ran\]ds.$$
\ep

We now consider the case that $[A,B]$ is not assumed to be stabilizable and the following are not assumed:
\bel{b-q-rho}b(\cd),q(\cd)\in L^1(0,\i;\dbR^n),\q\rho(\cd)\in L^2(0,\i;\dbR^m).\ee
Then \rf{U non-empty} could fail. For such cases, there are two possible approaches to deal with such a situation. Let us briefly look at them.

\ms

\it Approach 1. Ces\`aro mean. \rm Define
\bel{J}\h J_{_T}(x;u(\cd))={1\over T}\int_0^Tg(s,X(s),u(s))ds,\qq T>0.\ee
Then for each $T>0$, one poses the following problem:

\ms

\bf Problem (LQ)$^0_{_T}$. \rm For each $x\in\dbR^n$, find a $\bar u_{_T}(\cd)\in\sU_0[0,T]\equiv L^2(0,T;\dbR^m)$ such that
$$\h J_{_T}(x;\bar u_{T}(\cd))=\inf_{u(\cd)\in\sU_0[0,T]}J_{_T}(x;u(\cd))\equiv\h V_{_T}(x).$$

\ms

It is a hope that as $T\to\i$, the value function $\h V_{_T}(\cd)$ as well as the optimal control $\bar u_{_T}(\cd)$ (assume it exists) are convergent. The following simple illustrative example tells us that the above approach could fail.

\bex{ex-2.7} \rm Consider controlled system:
$$\dot X(s)=AX(s)+Bu(s),$$
with
$$A=\begin{pmatrix}1&0\\ 1&2\end{pmatrix},\qq B=\begin{pmatrix}1\\-1\end{pmatrix}.$$
Note that
$$AB=\begin{pmatrix}1&0\\ 1&2\end{pmatrix}\begin{pmatrix}1\\-1\end{pmatrix}=\begin{pmatrix}1\\-1\end{pmatrix}=B.$$
Thus, $[A,B]$ is not controllable. We have
$$\ba{ll}
\ns\ds X(s)=e^{As}x+\int_0^te^{A(s-\t)}Bu(\t)d\t=\begin{pmatrix}e^s&0\\ e^{2s}-e^s&e^{2s}\end{pmatrix}x+\int_0^s\begin{pmatrix}e^{s-\t}&0\\ e^{2(s-\t)}-e^{s-\t}&e^{2(s-\t)}\end{pmatrix}\begin{pmatrix}1\\-1\end{pmatrix}u(\t)
d\t\\
\ns\ds=\begin{pmatrix}e^sx_1+\int_0^se^{s-\t}u(\t)d\t\\
(e^{2s}-e^s)x_1+e^{2s}x_2-\int_0^se^{s-\t}u(\t)d\t\end{pmatrix}\equiv
\begin{pmatrix}X_1(s)\\
e^{2s}(x_1+x_2)-X_1(s)\end{pmatrix}.\ea$$
The running cost rate is given by the following:
$$g(X(s),u(s))=|X(s)|^2+|u(s)|^2.$$
For any $T>0$, by $2ab\les{1\over2}a^2+2b^2$ with $a=e^{2s}(x_1+x_2)$, $b=X_1(s)$, one has
$$\ba{ll}
\ns\ds\h J_{_T}(0,x;u(\cd))={1\over T}\int_0^T\(|X_1(s)|^2+[e^{2s}(x_1+x_2)-X_1(s)]^2
+|u(s)|^2\)ds\\
\ns\ds\qq\qq\q\ \ \ges{1\over T}\int_0^T\(\frac{1}{2}e^{4s}|x_1+x_2|^2+|u(s)|^2\)ds\ges{1\over8T}(e^{4T}-1)
|x_1+x_2|^2.\ea$$
Therefore, as long as $x_1\ne -x_2$, the limit of $\h V_{_T}(\cd)$ does not exist.

\ex

\it Approach 2. Abel mean. \rm For any $\l>0$, define the following discounted cost functional
$$J^\l(x;u(\cd))=\int_0^\infty e^{-\l s}g(s,X(s),u(s))ds.$$
Note that for the zero control $u_0(\cd)=0$, if we let $\l_A=\max\Re\si(A)>0$, then for any $\e>0$,
$$\sup_{s\ges0}\(e^{-(\l_A+\e)s}\|X(s;x,u_0(\cd))\|\)<\infty.$$
Hence, for $\l>2\l_A$, one has $J^\l(x;u_0(\cd))<\i$,
which implies that
\bel{}\sU_{ad}^\l[0,\infty)\equiv\Big\{u(\cd)\in\sU_{loc}[0,\infty)\bigm|J^\l(x;u(\cd))\hb{ exists }\Big\}\ne\varnothing.\ee
Now, one may take $\l$ small and still keep $\sU_{ad}^\l[0,\infty)\ne\varnothing$. For example, if $[A,B]$ is stabilizable, then $\sU_{ad}^0[0,\infty)\ne\varnothing$. In any case, for $\l>0$ large enough, we could pose the following problem.

\ms

\bf Problem (LQ)$^\l$. \rm For each $x\in\dbR^n$, find a $\bar u^\l(\cd)\in\sU_{ad}^\l[0,\infty)$ such that
$$J^\l(x;\bar u^\l(\cd))=\inf_{u(\cd)\in\sU_{ad}^\l[0,\infty)}J^\l(x;u(\cd))\equiv \h V^\l(x).$$

Again, it is a hope that as $\l\to0$, the value function $\h V^\l(\cd)$ as well as the optimal control $\bar u^\l(\cd)$ (if it exists) are convergent. However, let us look at Example \ref{ex-2.7} again. One has
$$e^{-\l s}\(|X(s)|^2+|u(s)|^2\)\ges e^{-\l s}\(\frac{1}{2}e^{4s}|x_1+x_2|^2+|u(s)|^2\)\ges \frac{1}{2}e^{(4-\l)s}|x_1+x_2|^2.$$
Hence, as long as $x_1\ne- x_2$, even the optimal value $V^\l(x)$ is not finite (if $0<\l<4$).

\ms

Note that in \cite{Buckdahn-Li-Quincampoix-Renault 2020}, both Ces\`aro and Abel means were considered for a more general class of nonlinear stochastic problems and the convergence of $\h V_{_T}(\cd)$ and $\h V^\l(\cd)$ were obtained, under the condition that the state $X(s)\equiv X(s;x,u(\cd))$ stays in a bounded set (depending on $x$). Our example does not satisfy such a condition. Also, in \cite{Basco-Cannarsa-Frankowska 2018}, for a more general nonlinear problem, the value function is allowed to be infinite and the theory was established on the set on which the value function is finite. However, for our example, the value function (without discount) is finite only on the set $x_1=-x_2$, and on which the value function is equal to zero. Hence, the results of \cite{Basco-Cannarsa-Frankowska 2018} are useless for such an example. From the above observations, we see that for the case that $J(x;u(\cd))$ is not well-defined (over infinite time horizons), we might want to take some other (possibly better) approaches. This naturally leads to the following definition for Problem (LQ) (see \cite{Carlson-Haurie-Leizarowitz 1991}).

\bde{overtaking} \rm For the control system \rf{state1}, let
\bel{J_T}J_{_T}(t,x;u(\cd))=\int_t^Tg(s,X(s),u(s))ds,\qq(t,x)\in[0,\i)\times\dbR^n,~u(\cd)
\in\sU[t,\i),\ee
where $\sU[t,\i)$ is a non-empty subset of $\sU_{loc}[t,\i)$.

\ms

(i) Control $u^*(\cd)\in\sU[t,\i)$ is called a {\it strong optimal control} of Problem (LQ) at the initial pair $(t,x)\in[0,\i)\times\dbR^n$ over $\sU[t,\i)$ if $J(t,x;u^*(\cd))$ is finite and
\bel{J2}J(t,x;u^*(\cd))\les J(t,x;u(\cd)),\qq\forall u(\cd)\in\sU[t,\i).\ee

(ii) Control $u^*(\cd)\in\sU[t,\i)$ is called an {\it overtaking optimal control} of Problem (LQ) at the initial pair $(t,x)\in[0,\i)\times\dbR^n$ over $\sU[t,\i)$ if
\bel{limsup<0}\limsup_{T\to\infty}\[J_{_T}(t,x;u^*(\cd))-J_{_T}(t,x;u(\cd))\]\les0,\qq\forall\ u(\cd)\in\sU[t,\infty).\ee

(iii) Control $u^*(\cd)\in\sU[t,\infty)$ is called a {\it weakly overtaking optimal control} of Problem (LQ) at the initial pair $(t,x)\in[0,\i)\times\dbR^n$ over $\sU[t,\i)$ if
\bel{liminf<0}\liminf_{T\to\infty}\[J_{_T}(t,x;u^*(\cd))-J_{_T}(t,x;u(\cd))\]\les0,\qq\forall\ u(\cd)\in\sU[t,\infty).\ee

\ede

Note that for the notions of (weakly) overtaking optimal control, the functional $J(t,x;u^*(\cd))$ is not required to be well-defined, since we only need $J_{_T}(t,x;u^*(\cd))$ to be finite. It is clear that in the above definition, we have the following implications: (i) $\Ra$  (ii) $\Ra$ (iii). The overtaking optimal control problem is by no means trivial; both the existence and non-existence of overtaking optimal controls are very subtle. Note that in the above definition, $\sU[t,\i)$ is any non-empty subset of $\sU_{loc}[t,\i)$ which does not have to be a linear space.

\section{An Ideal Case}

In this section, we are going to present an ideal case for which one could obtain the existence of overtaking optimal controls. First, we present the following simple result which shows an interesting power of overtaking optimality.

\bp{ignore} \sl Let Problem (LQ) admit an open-loop overtaking optimal control $u^*(\cd)\in\sU[t,\i)$ at initial pair $(t,x)\in[0,\i)\times\dbR^n$. Let $\f(\cd)\in L^1_{loc}(t,\i;\dbR)$. Define
$$J_{_T}^\f(t,x;u(\cd))=J_{_T}(t,x;u(\cd))+\int_t^T\f(s)ds\equiv\int_t^T[g(s,X(s),u(s))+\f(s)]ds.$$
Then the LQ problem with the running cost rate function $g(s,X,u)+\f(s)$ admits the same overtaking optimal control.

\ep

\it Proof. \rm If $u^*(\cd)\in\sU[t,\infty)$ is an overtaking optimal control of the original problem, then for any $u(\cd)\in\sU[t,\i)$, it holds
$$\ba{ll}
\ns\ds\limsup_{T\to\i}\[J_{_T}^\f(x;u^*(\cd))-J_{_T}^\f(x;u(\cd))\]
=\limsup_{T\to\i}\[J_{_T}(x;u^*(\cd))-J_{_T}(x;u(\cd))\]\les0.\ea$$
This proves our conclusion. \endpf

\ms

It is clear that a similar result holds true for weakly overtaking optimal controls. According to the above result, we see that one can drop or add any (locally integrable) terms independent of the state and control in the running cost rate function without changing the existence/non-existence of (weakly) overtaking optimal controls.

\ms

Since we do not assume either controllability or stabilizability of system $[A,B]$, the decomposition of the system mentioned in the introduction is very useful. Let us now make it more precise. Let $\dbH_0=\span\big\{\sR(A^kB)\bigm|k\ges0\big\}$ which is the controllability space of system $[A,B]$. We assume that $\dim\dbH_0=\ell$. Let $\Pi:\dbR^n\to\dbH_0$ be the orthogonal projection. Then $\Pi^\perp=I-\Pi:\dbR^n\to(\dbH_0)^\perp$ is also an orthogonal projection. For any initial pair $(t,x)\in[0,\i)\times\dbR^n$, and $u(\cd)\in\sU_{loc}[t,\i)$, let $X(\cd)$ be the corresponding state process. Denote
\bel{X_Pi}X_{_\Pi}(s)=\Pi X(s),\q X_{_{\Pi^\perp}}(s)=\Pi^\perp X(s),\qq s\in[t,\i).\ee
Then
\bel{dot X_Pi}\dot X_{_\Pi}(s)=A_{_\Pi}X_{_\Pi}(s)+B_{_\Pi}u(s)+b_{_\Pi}(s),\ee
where
\bel{A_Pi}A_{_\Pi}=\Pi A\Pi,\qq B_{_\Pi}=\Pi B,\qq b_{_\Pi}(s)=\Pi A\Pi^\perp X_{_{\Pi^\perp}}(s)+\Pi b(s).\ee
Clearly, system $[A_{_\Pi},B_{_\Pi}]$ is controllable on $\dbH_0$. On the other hand, it is ready to see that
\bel{R(A)}\sR(A\Pi),\sR(B)\subseteq\sR(\Pi)=\dbH_0\qq\Ra\qq\Pi^\perp A\Pi=0,~\Pi^\perp B=0.\ee
Thus,
\bel{dot X_Pi*}\dot X_{_{\Pi^\perp}}(s)=A_{_{\Pi^\perp}}X_{_{\Pi^\perp}}(s)+b_{_{\Pi^\perp}}(s),\ee
where
\bel{A_Pi*}A_{_{\Pi^\perp}}=\Pi^\perp A\Pi^\perp,\qq b_{_{\Pi^\perp}}(s)=\Pi^\perp b(s).\ee
From this, we see that $X_{_{\Pi^\perp}}(\cd)$ is not affected by the control
$u(\cd)$. By the controllability of $[A_{_\Pi},B_{_\Pi}]$, we may take a stabilizer $\Th:\dbH_0\to\dbR^m$ of this subsystem such that $\si\big(A_{_\Pi}+B_{_\Pi}\Th\big)\subseteq\dbC^-$. By taking the control of the following form:
\bel{feedback}u(s)=\Th X_{_\Pi}(s)+v(s),\qq s\ges t,\ee
and combining \rf{dot X_Pi} and \rf{feedback}, we have
\bel{dot XX}\dot X_{_\Pi}(s)=(A_{_\Pi}+B_{_\Pi}\Th)X_{_\Pi}(s)+B_{_\Pi}v(s)+b_{_\Pi}(s).\ee
We now look at the running cost rate:
\bel{g*}\ba{ll}
\ns\ds g(s,X(s),u(s))=\lan QX(s),X(s)\ran+2\lan SX(s),u(s)\ran+\lan Ru(s),u(s)\ran+2\lan q(s),X(s)\ran+2\lan\rho(s),u(s)\ran\\
\ns\ds=\lan Q[X_{_\Pi}(s)+X_{_{\Pi^\perp}}(s)],X_{_\Pi}(s)+X_{_{\Pi^\perp}}(s)\ran+2\lan S[X_{_\Pi}(s)+X_{_{\Pi^\perp}}(s)],\Th X_{_\Pi}(s)+v(s)\ran\\
\ns\ds\qq+\lan R[\Th X_{_\Pi}(s)+v(s)],\Th X_{_\Pi}(s)+v(s)\ran+2\lan q(s),X_{_\Pi}(s)+X{_{\Pi^\perp}}(s)\ran+2\lan\rho(s),\Th X_{_\Pi}(s)+v(s)\ran\\
\ns\ds=\lan Q_{_\Pi}^\Th X_{_\Pi}(s),X_{_\Pi}(s)\ran\1n+\1n2\lan S_{_\Pi}^\Th X_{_\Pi}(s),v(s)\ran\1n+\1n\lan Rv(s),v(s)\ran+2\lan q_{_\Pi}^\Th(s),X_{_\Pi}(s)\ran\1n+\1n2\lan\rho_{_\Pi}(s),v(s)\ran
+\f_{_\Pi}(s)\\
\ns\ds\equiv g_{_\Pi}^\Th(s,X_{_\Pi}(s),v(s))+\f_{_\Pi}(s),\ea\ee
where
$$\ba{ll}
\ns\ds Q_{_\Pi}^\Th=\Pi(Q+S^\top\Th+\Th^\top S+\Th^\top R\Th)\Pi,\qq S_{_\Pi}^\Th=S\Pi+R\Th,\\
\ns\ds q_{_\Pi}^\Th(s)=\Pi q(s)+\Pi
(Q+\Th^\top S)\Pi^\perp X_{_{\Pi^\perp}}(s)+\Th^\top\rho(s),\qq\rho_{_\Pi}(s)=\rho(s)+S\Pi^\perp
X_{_{\Pi^\perp}}(s),\\
\ns\ds\f_{_\Pi}(s)=\lan QX_{_{\Pi^\perp}}(s),X_{_{\Pi^\perp}}(s)\ran+2\lan q(s),X_{_{\Pi^\perp}}(s)\ran.\ea$$

By Proposition \ref{ignore}, it is seen that if we consider the LQ problem with the state equation \rf{dot XX} and running cost rate function $g_{_\Pi}^\Th(s,X_{_\Pi}(s),v(s))$, call it Problem (LQ)$_\Pi^\Th$, then it admits a (weakly) overtaking optimal control if and only if so does the original Problem (LQ). Note that we do not claim the equivalence of strong optimality between Problems (LQ) and (LQ)$^\Th_\Pi$ as that requires the original cost functional $J(t,x;u^*(\cd))$ to be finite, which we might not have.

\ms

We now look at the most ideal case of Problem (LQ) for which it admits overtaking
optimal controls. Consider an LQ problem with the state equation
\bel{state Pi}\left\{\2n\ba{ll}
\ds\dot X_{_\Pi}(s)=A_{_\Pi}X_{_\Pi}(s)+B_{_\Pi}u(s)+b_{_\Pi}(s),\qq s\ges t,\\
\ns\ds X_{_\Pi}(t)=\Pi x\equiv x_{_\Pi},\ea\right.\ee
and with the running cost rate function ($\Th=0$ and the superscript $\Th$ is omitted below)
\bel{g_Pi}g_{_\Pi}(s,X_{_\Pi},u)=\lan Q_{_\Pi}X_{_\Pi},X_{_\Pi}\ran+2\lan S_{_\Pi}X_{_\Pi},u\ran+\lan Ru,u\ran+2\lan q_{_\Pi}(s),X_{_\Pi}\ran+2\lan\rho_{_\Pi}(s),u\ran.\ee
We let
\bel{J^Pi}J^\Pi(t,x_{_\Pi};u(\cd))=\int_t^\i g_{_\Pi}(s,X_{_\Pi}(s),u(s))ds,\qq
\forall(t,x_\Pi)\in[0,\i)\times\dbH_0,\ee
and let
$$\sU_{_\Pi}[t,\i)=\Big\{u(\cd)\in\sU_{loc}[t,\i)\bigm|J^\Pi(t,x_{_\Pi};u(\cd))\hb{ is
well-defined}\Big\}.$$
Then we may pose the following LQ problem.

\ms

\bf Problem (LQ)$_{_\Pi}$. \rm For any initial pair $(t,x_{_\Pi})\in[0,\i)\times\dbH_0$,
find a $\bar u_{_\Pi}(\cd)\in\sU_{_\Pi}[t,\i)$ such that
$$J^\Pi(t,x_{_\Pi};\bar u_{_\Pi}(\cd))=\inf_{u(\cd)\in\sU_{_\Pi}[t,\i)}J^\Pi(t,x_{_\Pi};u(\cd)).$$

\ms

We have the following result.

\bt{sufficient2} \sl Suppose
\bel{bqr Pi}b_{_\Pi}(\cd),q_{_\Pi}(\cd)\in L^1(0,\i;\dbH_0),\qq\rho_{_\Pi}(\cd)\in
L^2(0,\i;\dbR^m).\ee
Suppose $u(\cd)\mapsto J^\Pi(\cd,\cd;u(\cd))$ is uniformly convex on
$\sU_{_\Pi}[t,\i)$, which is true if \rf{QSR} holds. Then
\bel{U=U}\sU_{_\Pi}[t,\i)=\sU_0[t,\i),\qq\forall t\ges0,\ee
and for any initial pair $(t,x)\in[0,\i)\times\dbR^n$, by letting $x_{_\Pi}=\Pi x$, Problem (LQ)$_{_\Pi}$ admits a unique optimal control $\bar u_{_\Pi}(\cd)\in\sU_0[t,\i)$ at $(t,x_{_\Pi})$. This $\bar u_{_\Pi}(\cd)$ must be an overtaking optimal control of Problem (LQ) at $(t,x)$ over $\sU_0[t,\i)$. Furthermore, $\bar u_{_\Pi}(\cd)$ can be obtained through the solution to an algebraic Riccati equation and an ODE on $[t,\i)$.

\et

\it Proof. \rm First of all, by the controllability of $[A_{_\Pi},B_{_\Pi}]$, it is not hard to show \rf{U=U}. Next, according to Proposition \ref{LQ results}, for any $(t,x)\in[0,\i)\times\dbR^n$ with $x_{_\Pi}=\Pi x$, Problem (LQ)$_\Pi$ admits a unique optimal control $\bar u_{_\Pi}(\cd)\in\sU_0[t,\i)$. We let $\bar X(\cd)=X(\cd\,;t,x,\bar u_{_\Pi}(\cd))$ and $\bar X_{_\Pi}(\cd)=X_{_\Pi}(\cd\,;t,x_{_\Pi},\bar u_{_\Pi}(\cd))$.
Then by the optimality of $\bar u_{_\Pi}(\cd)$, we have that for any $u(\cd)
\in\sU_{_\Pi}[t,\i)$,
$$\ba{ll}
\ns\ds J_T(t,x;\bar u_{_\Pi}(\cd))=\int_t^Tg(s,\bar X(s),\bar u_{_\Pi}(s))ds
=\int_t^T\(g_{_\Pi}(s,\bar X_{_\Pi}(s),\bar u_{_\Pi}(s))+\f_{_\Pi}(s)\)ds\\
\ns\ds\qq\qq\qq\ =J^\Pi(t,x_{_\Pi};\bar u_{_\Pi}(\cd))-\int_T^\i g_{_\Pi}(s,\bar X_{_\Pi}(s),\bar u_{_\Pi}(s))ds+\int_t^T\f_{_\Pi}(s)ds\\
\ns\ds\qq\qq\qq\ \les J^\Pi(t,x_{_\Pi};u(\cd))-\int_T^\i g_{_\Pi}(s,\bar X_{_\Pi}(s),\bar u_{_\Pi}(s))ds+\int_t^T\f_{_\Pi}(s)ds\\
\ns\ds\qq\qq\qq\ =J_T(t,x;u(\cd))+\int_T^\i g_{_\Pi}(s,X_{_\Pi}(s),u(s))ds-\int_T^\i g_{_\Pi}(s,\bar X_{_\Pi}(s),\bar u_{_\Pi}(s))ds.\ea$$
Note that since $u(\cd),\bar u_{_\Pi}(\cd)\in\sU_0[t,\i)$, we have $g_{_\Pi}(\cd\,,X_{_\Pi}(\cd),u(\cd))$ and $g_{_\Pi}(\cd\,,\bar X_{_\Pi}(\cd),\bar u_{_\Pi}(\cd))$ are integrable on $[t,\i)$. Thus, the last two terms on the right-hand side go to 0 as $T\to\i$. As a result, one has
$$\limsup_{T\to\i}\[J_{_T}(t,x;\bar u_{_\Pi}(\cd))-J_{_T}(t,x;u(\cd))\]\les0.$$
This means that $\bar u_{_\Pi}(\cd)$ is an overtaking optimal control in $\sU_{_\Pi}[t,\i)$. \endpf

\ms

We note that condition \rf{bqr Pi} gives some kind of compatibility among $b(\cd),q(\cd),\rho(\cd)$ and $X_{_{\Pi^\perp}}(\cd)$. In particular, if
\bel{3.14}\Pi Q\Pi^\perp=0,\qq S\Pi^\perp=0,\qq\Pi b(\cd),\Pi q(\cd)\in L^1(0,\i;\dbH_0),\qq\rho(\cd)\in L^2(0,\i;\dbR^m),\ee
then \rf{bqr Pi} holds. This amounts to saying that the LQ problem restricted on the space $\dbH_0$ is solvable and the restrictions of the system as well as the cost functional on $\dbH_0$ and $\dbH_0^\perp$ are decoupled. We emphasize here that no controllability condition is assumed for the original system $[A,B]$.

\ms

In the rest of this section, let us look at Example \ref{ex-2.7} again. As we already pointed out, $[A,B]$ is not controllable, therefore, we may let
$$\dbH_0=\span\big\{\sR(B),\sR(AB)\big\}=\Big\{\l\begin{pmatrix}1\\-1\end{pmatrix}\bigm|
\l\in\dbR\Big\}\equiv\span\Big\{{\Be_1-\Be_2\over\sqrt2}\Big\},$$
where $\Be_1=(1,0)^\top$, $\Be_2=(0,1)^\top$. Clearly,
$$\dbH_0^\perp=\Big\{\l\begin{pmatrix}1\\ 1\end{pmatrix}\bigm|\l\in\dbR\Big\}=\span\{
{\Be_1+\Be_2\over\sqrt2}\}.$$
In this case, we have
$$x_{_\Pi}\equiv\Pi x=\lan x,{\Be_1-\Be_2\over\sqrt2}\ran{\Be_1-\Be_2\over\sqrt2}=(x_1-x_2){\Be_1-\Be_2\over2}
=\begin{pmatrix}{1\over2}&-{1\over2}\\-{1\over2}&{1\over2}\end{pmatrix}x,
\qq\forall x\in\dbR^2,$$
and
$$x_{_{\Pi^\perp}}\equiv\Pi^\perp x=\lan x,{\Be_1+\Be_2\over\sqrt2}\ran{\Be_1+\Be_2\over\sqrt2}=(x_1+x_2){\Be_1+\Be_2\over2}
=\begin{pmatrix}{1\over2}&{1\over2}\\{1\over2}&{1\over2}\end{pmatrix}x,
\qq\forall x\in\dbR^2.$$
Then
$$\dot X_{_\Pi}(s)=\begin{pmatrix}{1\over2}&-{1\over2}\\-{1\over2}&{1\over2}\end{pmatrix}
\begin{pmatrix}1&0\\1&2
\end{pmatrix}\begin{pmatrix}{1\over2}&-{1\over2}\\-{1\over2}&{1\over2}\end{pmatrix}
X_{_\Pi}(s)+\begin{pmatrix}{1\over2}&-{1\over2}\\-{1\over2}&{1\over2}\end{pmatrix}
\begin{pmatrix}1\\-1\end{pmatrix}u(s)=X_{_\Pi}(s)+\begin{pmatrix}1\\-1\end{pmatrix}u(s),$$
and
$$\dot X_{_{\Pi^\perp}}(s)=\begin{pmatrix}{1\over2}&{1\over2}\\{1\over2}&{1\over2}\end{pmatrix}
\begin{pmatrix}1&0\\1&2\end{pmatrix}
\begin{pmatrix}{1\over2}&{1\over2}\\{1\over2}&{1\over2}\end{pmatrix}
X_{_{\Pi^\perp}}(s)=2X_{_{\Pi^\perp}}(s).$$
The running cost rate is
$$g(s,X(s),u(s))=|X_{_\Pi}(s)+X_{_{\Pi^\perp}}(s)|^2+|u(s)|^2=|X_{_\Pi}(s)|^2
+|X_{_{\Pi^\perp}}(s)|^2+|u(s)|^2.$$
Hence, the overtaking optimal control can be obtained by solving the following one-dimensional LQ problem (in $\dbH_0$) whose state equation is: (denoting
$\Be={\Be_1-\Be_2\over\sqrt{2}}$)
$$\dot X_{_\Pi}(s)=X_{_\Pi}(s)+\sqrt2\,u(s)\Be,$$
with the cost functional
$$J^\Pi(t,x_{_\Pi};u(\cd))=\int_t^\i\(|X_{_\Pi}(s)|^2+|u(s)|^2\)ds.$$
More precisely, the corresponding Riccati equation reads
$$2P+1-2P^2=0,$$
whose positive solution is given by $P={1+\sqrt3\over2}$. Thus,
$$u(s)=-{1+\sqrt3\over\sqrt2}\lan X_{_\Pi}(s),\Be\ran,\qq s\in[t,\i).$$
Solving the closed loop system, we can obtain
$$\bar u(s)=e^{-\sqrt3 s}\lan x_{_\Pi},\Be\ran=\frac{x_1-x_2}{\sqrt{2}}e^{-\sqrt3 s},\qq s\in[t,\i).$$
According to the above result, this $\bar u(\cd)$ is an overtaking optimal control for the corresponding LQ problem at $(t,x)$ over $\sU_0[t,\i)$.

\section{Some General Considerations}

We now would like to consider general cases for which compatibility conditions \rf{bqr Pi} (or something like \rf{3.14}) are not assumed. To begin, let us present the following result showing
that it is not very restrictive by imposing some additional conditions for the
running cost rate function.

\bp{bounded below} \sl Let Problem (LQ) admit an overtaking optimal control $u^*(\cd)\in\sU[t,\i)$ at some initial pair $(t,x)\in[0,\i)\times\dbR^n$ over $\sU_0[t,\i)$. Then, for any $\h T>t$, there exists an interval $[S,T]$ with $\h T\les S<T<\i$ such that the LQ problem posed on $[S,T]$, called Problem (LQ)$_{[S,T]}$, must be finite, i.e.,
$$\inf_{u(\cd)\in\sU[S,T]}J_{_T}(S,x;u(\cd))\equiv\inf_{u(\cd)\in\sU[S,T]}\int_S^T
g(s,X(s),u(s))ds>-\i,\qq\forall(t,x)\in[S,T]\times\dbR^n.$$
Consequently, the map $u(\cd)\mapsto J_{_T}(S,x;u(\cd))$ must be convex.

\ep

\it Proof. \rm Suppose otherwise. Then there exists a strictly increasing sequence $T_k\to\i$ such that Problem (LQ)$_{[T_k,T_{k+1}]}$ is not finite, i.e.,
\bel{non-finite*}\inf_{u(\cd)\in\sU[T_k,T_{k+1}]}J_{_{T_{k+1}}}(T_k,x;u(\cd))=-\i.\ee
Now, suppose Problem (LQ) admits an overtaking optimal pair $(X^*(\cd),
u^*(\cd))$. Then the following holds:
$$\limsup_{T\to\infty}\[J_{_T}(x;u^*(\cd))-J_{_T}(x;u(\cd))\]\les0,\qq\forall u(\cd)\in\sU[t,\i).$$
By \rf{non-finite*}, we can find a $u_0(\cd)\in\sU[T_0,T_1]$ such that
$$J_{_{T_1}}(T_0,x;u_0(\cd))<J_{_{T_1}}(T_0,x;u^*(\cd))-1.$$
Next, we find $u_1(\cd)\in\sU[T_1,T_2]$ such that
$$J_{_{T_2}}(T_1,X(T_1);u_1(\cd))<J_{_{T_2}}(T_1,X^*(T_1);u^*(\cd))-1.$$
By induction, we can find $u_k(\cd)\in\sU[T_k,T_{k+1}]$ such that
$$J_{_{T_{k+1}}}(T_k,X(T_k);u_k(\cd))<J_{_{T_{k+1}}}(T_k,X^*(T_k);u^*(\cd))-1.$$
We patch $u_k(\cd)$ together to get a $u(\cd)\in\sU[t,\i)$ and we have
$$J_{_{T_k}}(x;u(\cd))=\sum_{i=1}^kJ_{_{T_i}}(T_{i-1},X(T_{i-1});u_{i-1}(\cd))<
\sum_{i=1}^kJ_{_{T_i}}(T_{i-1},X^*(T_{i-1});u^*(\cd))-k=J_{_{T_k}}(x;u^*(\cd))-k.$$
This leads to
$$\limsup_{T\to\i}\[J_{_T}(x;u^*(\cd))-J_{_T}(x;u(\cd))\]\ges\lim_{k\to\infty}
\[J_{_{T_k}}(x;u^*(\cd))-J_{_{T_k}}(x;u(\cd))\]=\i,$$
a contradiction. Finally, by the finiteness of the LQ problem on $[S,T]$ and the assumption that the quadratic weighting matrices of the running cost rate and the coefficients of the system are constants, we must have the convexity of $u(\cd)\mapsto J_{_T}(S,x;u(\cd))$ (see \cite{Sun-Yong 2020}). \endpf

\ms

The above result tells us that for overtaking optimal control problems, it is not extremely restrictive to assume that for any $T>0$, the map $u(\cd)\mapsto J_{_T}(t,x;u(\cd))$ is convex on $\sU[t,T]$, or even assume \rf{QSR}. It is known that under the above conditions, on any $[S,T]$, the corresponding LQ problem is closed-loop solvable (\cite{Sun-Yong 2020}).

\ms

We now consider state equation \rf{state1} and running cost rate \rf{cost rate}, with $[A,B]$ being controllable, no global integrability condition is assumed for $b(\cd),q(\cd),\rho(\cd)$. In what follows, \rf{QSR} will be assumed. We now would like to make a further reduction. Let us first look at the following (suppressing $s$ in $X(s)$ and $u(s)$)
$$\ba{ll}
\ns\ds g(s,X,u)=\lan QX,X\ran+2\lan SX,u\ran+\lan Ru,u\ran+2\lan q,X\ran+2\lan\rho,u\ran\\
%
%
\ns\ds\qq\qq\ =\lan QX,X\ran+2\lan q,X\ran+|R^{1\over2}u+R^{-{1\over2}}
(SX+\rho)|^2-|R^{-{1\over2}}(SX+\rho)|^2.\ea$$
Denote
$$\h u=R^{1\over2}[u+R^{-1}(SX+\rho)].$$
Then
$$u=R^{-{1\over2}}\h u-R^{-1}(SX+\rho),$$
and the state equation \rf{state1} becomes
$$\dot X=AX+B\big[R^{-{1\over2}}\h u-R^{-1}(SX+\rho)\big]+b=(A-BR^{-1}S)X+BR^{-{1\over2}}\h u+b-BR^{-1}\rho.$$
Next, we decompose the state $X(\cd)$ in the following manner:
$$X(\cd)=\wt X(\cd)+X_0(\cd),$$
where $X_0(\cd)$ is the solution to the following
\bel{X_0}\left\{\2n\ba{ll}
\ds\dot X_0(s)=(A-BR^{-1}S)X_0(s)+b(s)-BR^{-1}\rho(s),\qq s\ges t,\\
\ns\ds X_0(t)=0,\ea\right.\ee
which only depends on the non-homogenous term $b(\cd)-BR^{-1}\rho(\cd)$,
independent of the control $\h u(\cd)$ and the initial state $x$; and $\wt X(\cd)$ is the solution to the following equation:
$$\left\{\2n\ba{ll}
\ds\dot{\wt X}(s)=\big[A+BR^{-1}(\Th-S)\big]\wt X(s)
+BR^{-{1\over2}}\wt u(s),\qq s\ges t,\\
\ns\ds\wt X(t)=x,\ea\right.$$
where $\Th\in\dbR^{m\times n}$ is so chosen that $\si\big(A+BR^{-1}(\Th-S)\big)\subseteq\dbC^-$ and $\h u=R^{-{1\over2}}\Th\wt X+\wt u$. Then
$$\ba{ll}
\ns\ds g(s,X,u)=\lan Q(\wt X+X_0),\wt X+X_0\ran+2\lan q,\wt X+X_0\ran
+|R^{-{1\over2}}\Th\wt X+\wt u|^2-|R^{-{1\over2}}(S\wt X+SX_0+\rho)|^2\\
\ns\ds\qq\qq\ =\lan Q\wt X,\wt X\ran+2\lan QX_0,\wt X\ran+\lan QX_0,X_0\ran+2\lan q,
\wt X\ran+2\lan q,X_0\ran+\lan\Th^\top R^{-1}\Th\wt X,\wt X\ran
\\
\ns\ds\qq\qq\qq+2\lan R^{-{1\over2}}\Th\wt X,\wt u\ran+|\wt u|^2-\lan S^\top R^{-1}S\wt X,\wt X\ran-2\lan S^\top R^{-1}(SX_0+\rho),\wt X\ran-|R^{-{1\over2}}
(SX_0+\rho)|^2\\
\ns\ds\qq\qq\ \equiv\lan\wt Q\wt X(s),\wt X(s)\ran+2\lan\wt S\wt X(s),\wt u(s)\ran+|\wt u|^2
+2\lan\wt q(s),\wt X(s)\ran+\wt\f(s)\equiv\wt g(s,\wt X(s),\wt u(s))
+\wt\f(s),\ea$$
where
$$\ba{ll}
\ns\ds\wt Q=Q+\Th^\top R^{-1}\Th-S^\top R^{-1}S,\q\wt S=R^{-{1\over2}}\Th,\\
\ns\ds\wt q(s)=q(s)+(Q-S^\top R^{-1}S)X_0(s)
-S^\top R^{-1}\rho(s),\\
\ns\ds\wt\f(s)=\lan QX_0(s),X_0(s)\ran+2\lan q(s),X_0(s)\ran
-|R^{-{1\over2}}[SX_0(s)+\rho(s)]|^2.\ea$$
If we let
$$\wt A=A+BR^{-1}(\Th-S),\qq\wt B=BR^{-{1\over2}},$$
then the state equation becomes
$$\left\{\2n\ba{ll}
\ds\dot{\wt X}(s)=\wt A\wt X(s)+\wt B\wt u(s),\qq s\ges t,\\
\ns\ds\wt X(t)=x.\ea\right.$$
Since $\wt\f(\cd)$ is independent of $(\wt X(\cd),\wt u(\cd))$, by Proposition
\ref{ignore}, we may drop it, and take the running cost rate function
$\wt g(s,\wt X,\wt u)$. Note that under \rf{QSR},
$$\wt Q-\wt S^\top\wt S=Q-S^\top R^{-1}S\ges0.$$
The nonhomogeneous term $b(\cd)$ appears in $\wt q(\cd)$ through $X_0(\cd)$. Also, $X_0(\cd)$ depends not only on $b(\cd)$, but also on $\rho(\cd)$.

\ms

The above reduction tells us that without loss of generality, we may consider the state equation
\bel{state3}\left\{\2n\ba{ll}
\ds\dot X(s)=AX(s)+Bu(s),\qq s\in[t,\i),\\
\ns\ds X(t)=x,\ea\right.\ee
with the running cost rate
\bel{g}g(s,x,u)=\lan Qx,x\ran+2\lan Sx,u\ran+|u|^2+2\lan q(s),x\ran,\ee
where, we introduce the following hypothesis.

\ms

{\bf(H)} $[A,B]$ is controllable, $A$ is stable with
\bel{e^A}|e^{As}|\les Me^{-\m s},\qq\forall s\ges0,\ee
for some $M\ges1$ and $\m>0$, and
\bel{I-SS>0,q}Q-S^\top S\ges0,\qq q(\cd)\in L^1_{loc}(0,\i;\dbR^n)\setminus L^1(0,\i;\dbR^n).\ee

Note that under (H), we have $\dbH_0=\dbR^n$ and $[A_\Pi,B_\Pi]=[A,B]$. Due to the fact that $q(\cd)\notin L^1(0,\i;\dbR^n)$, \rf{bqr Pi} fails. In the rest of the paper, we concentrate on the problem associated with
state equation \rf{state3} and running cost rate function \rf{g}, which will
still be simply referred to as Problem (LQ).
For such a problem, we have the following result.

\bt{Th-4.2} \sl Let {\rm(H)} hold. Then for each $(t,x)\in[0,\i)\times\dbR^n$ and $u(\cd)\in\sU_0[t,\i)$,
\bel{F}\ba{ll}
\ns\ds F_0(s)x=Se^{A(s-t)}x+\int_s^\i B^\top e^{A^\top(\t-s)}Qe^{A(\t-t)}xd\t,\\
\ns\ds F_1[u(\cd)](s)=SX_0(s)+\2n
\int_s^\i\2n B^\top e^{A^\top(\t-s)}\big[QX_0(\t)+S^\top u(\t)\big]d\t,\ea\qq s\in[t,\i),\ee
exist satisfying
\bel{F*}\ba{ll}
\ns\ds\int_t^\i\big|F_0(s)x\big|^2ds\les
\(\|S\|M+{\|B\|M^2\|Q\|\over2\m}\)^2{|x|^2\over2\m},\\
\ns\ds\int_t^\i\big|F_1[u(\cd)](s)\big|^2ds\les{3\|B\|^2M^2\over\m^2}
\({\|B\|^2M^2\|Q\|^2\over\m^2}+2\|S\|^2\)\(\int_t^\i|u(s)|^2ds\),\ea\ee
where $X_0(\cd)=X_0(\cd\,;t,u(\cd))$ is given by the following:
\bel{}X_0(s)=\int_t^se^{A(s-\t)}Bu(\t)d\t,\qq s\in[t,\i).\ee
%
Further, let $\sU[t,\i)\subseteq\sU_0[t,\i)$ be convex. Then

\ms

{\rm(i)} $\bar u(\cd)\in\sU[t,\i)$ is an overtaking optimal control of Problem (LQ) at $(t,x)\in[0,\i)\times\dbR^n$ over $\sU[t,\i)$ if and only if
\bel{liminf**}0\les\liminf_{T\to\i}\int_t^T\lan\int_s^TB^\top e^{A^\top(\t-s)}q(\t)d\t+F_0(s)x+\bar u(s)+F_1[\bar u(\cd)](s),u(s)-\bar u(s)\ran ds,
\q\forall u(\cd)\in\sU[t,\i).\ee

\ms

{\rm(ii)} $\bar u(\cd)\in\sU[t,\i)$ is a weakly overtaking optimal control of Problem (LQ) at $(t,x)$ over $\sU[t,\i)$ if and only if
\bel{limsup**}0\les\limsup_{T\to\i}\int_t^T\lan\int_s^TB^\top e^{A^\top(\t-s)}q(\t)d\t+F_0(s)x+\bar u(s)+F_1[\bar u(\cd)](s),u(s)-\bar u(s)\ran ds,\q\forall u(\cd)\in\sU[t,\i).\ee

\et

\it Proof. \rm First, for any $u(\cd)\in\sU_0[t,\i)$, by \rf{e^A} and Young's inequality, we have
$$\ba{ll}
\ns\ds\int_t^\i\(\|S\|Me^{-\m s}+\|B\|M^2\|Q\|
\int_s^\i e^{-\m(\t-s)}e^{-\m(\t-t)}d\t\)^2ds|x|^2\\
\ns\ds=\int_t^\i\(\|S\|Me^{-\m s}+{\|B\|M^2\|Q\|\over2\m}e^{-\m(s-t)}\)^2ds|x|^2
=\(\|S\|M+{\|B\|M^2\|Q\|\over2\m}e^{\m t}\)^2{e^{-2\m t}\over2\m}|x|^2\\
\ns\ds=\(\|S\|Me^{-\m t}+{\|B\|M^2\|Q\|\over2\m}\)^2{|x|^2\over2\m}\les
\(\|S\|M+{\|B\|M^2\|Q\|\over2\m}\)^2{|x|^2\over2\m}.\ea$$
Next, due to
\bel{X_0}\ba{ll}
\ns\ds\int_t^\i|X_0(s)|^2ds=\int_t^\i\Big|\int_t^se^{A(s-\t)}B(\t)d\t\Big|^2ds\les M^2\|B\|^2\int_t^\i\(\int_t^se^{-\m(s-\t)}|u(\t)|d\t\)^2ds\\
\ns\ds\qq\qq\qq\ \ \les M^2\|B\|^2\(\int_t^\i e^{-\m\t}d\t\)^2\(\int_t^\i| u(\t)|^2d\t\)={M^2\|B\|^2\over\m^2}\(\int_t^\i|u(\t)|^2d\t\),\ea\ee
we have
$$\ba{ll}
\ns\ds\int_t^\i\Big|\int_s^\i B^\top e^{A^\top(\t-s)}[QX_0(\t)+S^\top u(\t)]d\t+SX_0(s)\Big|^2ds\\
\ns\ds\les\int_t^\i\(\int_s^\i\|B\|Me^{-\m(\t-s)}\big[\|Q\|\,|X_0(\t)|+\|S\|\,
|u(\t)|\big]d\t+\|S\|\,|X_0(s)|\)^2ds\\
%
%
%
\ns\ds\les3\|B\|^2M^2\(\1n\int_t^\i\2n e^{-\m s}ds\)^2\int_t^\i\(
\|Q\|^2|X_0(s)|^2\1n+\1n\|S\|^2|u(s)|^2\)ds+\1n3\|S\|^2\2n\int_t^\i\2n|X_0(s)|^2ds\1n\\
%
%
\ns\ds\les{3\|B\|^2M^2\over\m^2}\({\|B\|^2M^2\|Q\|^2\over\m^2}+2\|S\|^2\)
\(\int_t^\i|u(s)|^2ds\).\ea$$
As a result, $F_0(s)x$ and $F_1[u(\cd)](s)$ are will-defined and estimates \rf{F*} hold.

Now, suppose $\sU[t,\i)\subseteq\sU_0[t,\i)$ is convex. Then $\bar u(\cd)\in\sU[t,\i)$ is an overtaking optimal control of Problem (LQ) over $\sU[t,\i)$ if and only if for any $u(\cd)\in\sU[t,\i)$ and $\n\in[0,1]$,
\bel{limJ-J<0}\limsup_{T\to\i}\[J_{_T}\big(t,x;\bar u(\cd)\big)-J_{_T}\big(t,x;\bar u(\cd)+\n[u(\cd)-\bar u(\cd)]\big)\]\les0.\ee
Let us calculate the following
\bel{J-J*}\ba{ll}
\ns\ds J_{_T}\big(t,x;\bar u(\cd)\big)-J_{_T}\big(t,x;\bar u(\cd)+\n[u(\cd)
-\bar u(\cd)]\big)\\
\ns\ds=\int_t^T\(\lan Q\bar X(s),\bar X(s)\ran+2\lan S\bar X(s),\bar u(s)\ran+|\bar u(s)|^2+2\lan q(s),\bar X(s)\ran
-|Q[\bar X(s)+\n\xi(s)]|^2\\
\ns\ds\qq\q-2\lan S[\bar X(s)\1n+\1n\n\xi(s)],\bar u(s)\1n+\1n\n[u(s)\1n-\1n\bar u(s)]\ran\1n-\1n|\bar u(s)\1n+\1n\n[u(s)\1n-\1n\bar u(s)]|^2\1n-\1n2\lan q(s),\bar X(s)\1n+\1n\n\xi(s)\ran\)ds\\
%
%
%
\ns\ds=-2\n\int_t^T\lan q(s),\xi(s)\ran ds-2\n\int_t^T\(\lan Q\bar X(s)+S^\top\bar u(s),\xi(s)\ran+\lan S\bar X(s)+\bar u(s),u(s)-\bar u(s)\ran\)ds\\
\ns\ds\qq\q-\n^2\int_t^T\(|(Q-S^\top S)^{1\over2}\xi(s)|^2+|S\xi(s)+u(s)-\bar u(s)|^2\)ds,\ea\ee
where
$$\xi(s)=\int_t^se^{A(s-\t)}B[u(\t)-\bar u(\t)]d\t\equiv X_0(s;t,u(\cd)-\bar u(\cd)),\qq s\in[t,\i).$$
Note that
$$\int_t^T\2n\lan q(s),\xi(s)\ran ds\1n=\1n\int_t^T\3n\lan q(s),\2n\int_t^s\2n e^{A(s-\t)}B[u(\t)\1n-\1n\bar u(\t)]d\t\ran ds=\2n\int_t^T\2n\lan\int_s^T\2n B^\top e^{A^\top(\t-s)}q(\t)d\t,u(s)\1n-\1n\bar u(s)\ran ds,$$
and
$$\ba{ll}
\ns\ds\lim_{T\to\i}\int_t^T\(\lan Q\bar X(s)+S^\top\bar u(s),\xi(s)\ran+\lan S\bar X(s)+\bar u(s),u(s)-\bar u(s)\ran\)ds\\
%
%
\ns\ds=\int_t^\i\lan\int_s^\i B^\top
e^{A^\top(\t-s)}[Q\bar X(\t)+S^\top\bar u(\t)]d\t+S\bar X(s)+\bar u(s),u(s)-\bar u(s)\ran ds\\
\ns\ds\equiv\int_t^\i\lan F_0(s)x+\bar u(s)+F_1[\bar u(\cd)](s),u(s)-\bar u(s)\ran ds,\ea$$
where $\bar X_0(\cd)=X_0(\cd\,;t,\bar u(\cd))$.
Consequently, \rf{limJ-J<0} is equivalent to the following:
\bel{0<n}\ba{ll}
\ns\ds0\les\liminf_{T\to\i}\[2\n\int_t^T\lan\int_s^TB^\top e^{A^\top(\t-s)} q(\t)d\t+F_0(s)x+\bar u(s)+F_1[\bar u(\cd)](s),u(s)-\bar u(s)\ran ds\\
\ns\ds\qq\qq\q+\n^2\2n\int_t^\i\2n\(|(Q\1n-\1n S^\top S)^{1\over2}\xi(s)|^2\1n+\1n|S\xi(s)\1n+\1n u(s)\1n-\1n\bar u(s)|^2\)ds\],\q\forall u(\cd)\in\sU[t,\i),~\n\in[0,1].\ea\ee
Here, similar to \rf{J-J*}, we have the existence of the term following $\n^2$. Dividing $2\n$ and sending $\n\to0$, we see that \rf{liminf**} holds. Conversely, if \rf{liminf**} holds, then \rf{0<n} holds since the second term on the right-hand side is non-negative. Hence, $\bar u(\cd)$ is overtaking optimal.

\ms

(ii) The proof is the same, replacing $\liminf$ by $\limsup$ in the above. \endpf

\ms

We have the following simple and useful corollary.

\bc{Cor-4.3} \sl Let {\rm(H)} hold and $\sU[t,\i)\subseteq\sU_0[t,\i)$ be convex.

\ms

{\rm(i)} If $\bar u(\cd)\in\sU[t,\i)$ is an overtaking optimal control of Problem (LQ) at $(t,x)\in[0,\i)\times\dbR^n$ over $\sU[t,\i)$, then there exists a constant $C_0>0$ such that
\bel{liminf>}\liminf_{T\to\i}\int_t^T\lan\int_s^TB^\top e^{A^\top(\t-s)}q(\t)d\t,u(s)-\bar u(s)\ran ds\ges-C_0\(\int_t^\i|u(s)-\bar u(s)|^2ds\)^{1\over2},
\q\forall u(\cd)\in\sU[t,\i).\ee

\ms

{\rm(ii)} If $\bar u(\cd)\in\sU[t,\i)$ is a weakly overtaking optimal control of Problem (LQ) at $(t,x)\in[0,\i)\times\dbR^n$ over $\sU[t,\i)$, then there exists a constant $C_0>0$ such that
\bel{limsup>}\limsup_{T\to\i}\int_t^T\lan\int_s^TB^\top e^{A^\top(\t-s)}q(\t)d\t,u(s)-\bar u(s)\ran ds\ges-C_0\(\int_t^\i|u(s)-\bar u(s)|^2ds\)^{1\over2},
\q\forall u(\cd)\in\sU[t,\i).\ee

\ec

\it Proof. \rm Note that
\bel{C_0}\ba{ll}
\ns\ds\(\int_t^\i|F_0(s)x+\bar u(s)+F_1[\bar u(\cd)](s)|^2ds\)^{1\over2}\\
\ns\ds\les\(\|S\|M+{\|B\|M^2\|Q\|\over2\m}\){|x|\over\sqrt{2\m\,}}
+\[{\sqrt3\|B\|M\over\m}\({\|B\|^2M^2\|Q\|^2\over\m^2}
+2\|S\|^2\)^{1\over2}+1\]\(\int_t^\i|\bar u(s)|^2ds\)^{1\over2}.\ea\ee
By taking $C_0\equiv C_0(t,|x|,\bar u(\cd))$ to be the right-hand side of the above, we obtain our conclusions. \endpf

\section{Non-Existence of Overtaking Optimal Controls}

In this section, we are going to present some cases for which Problem (LQ) does not
admit overtaking optimal controls. Before stating and proving our results, let us
first make an observation. Suppose $q(\cd)$ grows exponentially, say, satisfying
\bel{e^a}0<q_0=\liminf_{s\to\i}e^{-\a s}|q(s)|\les\limsup_{s\to\i}e^{-\a s}|q(s)|=q_1<\i,\ee
for some $\a>0$. Then when $[A,B]$ is controllable, by making a state feedback first, we may make the system exponentially stable with any prespecified decay rate. Hence, for such a case, we may assume \rf{e^A} holds with $\m>\a$. Then
$$\int_s^TB^\top e^{A^\top(\t-s)}q(\t)d\t=B^\top e^{-A^\top s}\int_s^Te^{A^\top\t}q(\t)d\t\to B^\top e^{-A^\top s}\int_s^\i e^{A^\top\t}q(\t)d\t\equiv B^\top e^{-A^\top s}\h q(s),$$
as $T\to\i$, with
\bel{h q}\h q(s)=\int_s^\i e^{A^\top\t}q(\t)d\t,\qq s\in[0,\i).\ee
Also, since
$$e^{\m s}\int_s^\i e^{-\m\t}e^{\a\t}d\t={e^{\a s}\over\m-\a},$$
one expects that $B^\top e^{A^\top\cd}\h q(\cd)\notin L^1(0,\i;\dbR^m)$. On the other hand, it is possible that $q(\cd)$ grows faster than exponential functions, say, like $e^{s^2}$. In this case, one expects that
$$\int_0^\i|e^{A^\top\t}q(\t)|d\t=\i,$$
regardless of condition \rf{e^A}. Clearly, the above two cases are mutually exclusive. Now, let us present the following result.

\bt{Th-4.4} \sl Let {\rm(H)} hold and $q(\cd)\in L^1_{loc}(0,\i;\dbR^n)\setminus L^1(0,\i;\dbR^n)$ such that for some $\th_0\in S^m\equiv\{u\in\dbR^m\bigm||u|=1\}$,
\bel{th_0}\liminf_{T\to\infty}\big\lan\th_0,{\int_s^T B^\top e^{A^\top(\t-s)}q(\t)d\t
\over\big|\int_s^TB^\top e^{A^\top(\t-s)}q(\t)d\t\big|}\big\ran>0,\qq s\in[0,\i).\ee
Let $\sU[t,\i)\subseteq\sU_0[t,\i)$ be a convex set and $\bar u(\cd)$ be an interior point of $\sU[t,\i)$. Then $\bar u(\cd)$ must not be a weakly overtaking optimal control of Problem (LQ) if one of the following holds:

\ms

{\rm(i)} Let $e^{A^\top\1n\cd}q(\cd)\in L^1(0,\i;\dbR^n)$, and
\bel{q=infty}\limsup_{T\to\i}{1\over\sqrt T}\int_t^T\Big|\int_s^\i B^\top e^{A^\top(\t-s)}q(\t)d\t\Big|ds=\i.\ee

{\rm(ii)} It holds
\bel{q=infty*}\limsup_{T\to\i}\int_t^{T_0}\Big|\int_{T_0}^TB^\top e^{A^\top(\t-s)}q(\t)d\t\Big|ds=\i,\qq\forall T_0>t.\ee

\ms

\no Consequently, in one of the above cases, if $\sU[t,\i)\subseteq\sU_0[t,\i)$ is a subspace, then Problem (LQ) does not admit a weakly overtaking optimal control at $(t,x)$ over $\sU[t,\i)$.

\et

\it Proof. \rm Suppose $\bar u(\cd)$ is an interior point of $\sU[t,\i)$. Then, there exists a $\d>0$ such that
$$\bar u(\cd)+\eta(\cd)\in\sU[t,\i),\qq\forall\eta(\cd)\in\sU[t,\i),\hb{ with }\int_t^\i|\eta(s)|^2ds\les\d^2.$$
If $\bar u(\cd)$ is a weakly overtaking optimal control of Problem (LQ) over $\sU[t,\i)$, then
\bel{4.16}\limsup_{T\to\i}\int_t^T\lan\int_s^TB^\top e^{A^\top(\t-s)}q(\t)d\t,\eta(s)
\ran ds\ges-C\(\int_t^\i|\eta(s)|^2ds\)^{1\over2}\ges-C(1+\d)\equiv-C_0.\ee
Now, by \rf{th_0}, we can find a large $T_0>t$ such that
$$\big\lan\th_0,{\int_s^TB^\top e^{A^\top(\t-s)}q(\t)d\t
\over\big|\int_s^TB^\top e^{A^\top(\t-s)}q(\t)d\t\big|}\big\ran\ges\e>0,\qq T\ges T_0.$$
Then we take
$$\eta(s)=-{\d\over\sqrt{T_0-t\,}}\th_0{\bf1}_{[t,T_0)}(s),\qq s\in[t,\i).$$
Clearly,
$$\int_t^\i|\eta(s)|^2ds=\d^2,\qq\bar u(\cd)+\eta(\cd)\in\sU[t,\i).$$
In case (i), one has
$$\ba{ll}
\ns\ds\int_t^T\lan\int_s^TB^\top e^{A^\top(\t-s)}q(\t)d\t,\eta(s)
\ran ds=-{\d\over\sqrt{T_0-t\,}}\int_t^{T_0}\lan\int_s^TB^\top e^{A^\top(\t-s)}q(\t)d\t,\th_0\ran ds\\
%
%
\ns\ds\les-{\d\e\over\sqrt{T_0-t\,}}\[\int_t^{T_0}\Big|\int_s^\i B^\top e^{A^\top(\t-s)}q(\t)d\t\Big|ds-\int_t^{T_0}\Big|\int_T^\i B^\top e^{A^\top(\t-s)}q(\t)d\t\Big|ds\]\\
\ns\ds\les-{\d\e\over\sqrt{T_0-t\,}}\[\int_t^{T_0}\Big|\int_s^\i B^\top e^{A^\top(\t-s)}q(\t)d\t\Big|ds-\int_t^{T_0}\|Be^{-A^\top s}\|\int_T^\i|e^{A^\top\t}q(\t)|d\t ds\].\ea$$
Hence,
$$\limsup_{T\to\i}\int_t^T\lan\int_s^TB^\top e^{A^\top(\t-s)}q(\t)d\t,\eta(s)
\ran ds\les-{\d\e\over\sqrt{T_0-t\,}}\int_t^{T_0}\Big|\int_s^\i B^\top e^{A^\top(\t-s)}q(\t)d\t\Big|ds\to-\i,$$
as $T_0\to\i$, which contradicts \rf{4.16}.

\ms

In case (ii), one has
$$\ba{ll}
\ns\ds\int_t^T\lan\int_s^TB^\top e^{A^\top(\t-s)}q(\t)d\t,\eta(s)
\ran ds=-{\d\over\sqrt{T_0-t\,}}\int_t^{T_0}\lan\int_s^TB^\top e^{A^\top(\t-s)}q(\t)d\t,\th_0\ran ds\\
\ns\ds\les-{\d\over\sqrt{T_0-t\,}}\[\int_t^{T_0}\lan\int_s^{T_0}B^\top e^{A^\top(\t-s)}q(\t)d\t,\th_0\ran ds+\e\int_t^{T_0}\Big|\int_{T_0}^TB^\top e^{A^\top(\t-s)}q(\t)d\t\Big|ds\]\\
\ns\ds\les-{\d\over\sqrt{T_0-t\,}}\[\e\int_t^{T_0}\Big|\int_{T_0}^TB^\top e^{A^\top(\t-s)}q(\t)d\t\Big|ds-\int_t^{T_0}\Big|\int_s^{T_0}B^\top e^{A^\top(\t-s)}q(\t)d\t\Big|ds\].\ea$$
Hence, by \rf{q=infty}, we see that the first term on the right-hand side of the above is dominating and the second term stays finite as $T_0$ fixed. Consequently,
$$\limsup_{T\to\i}\int_t^T\lan\int_s^TB^\top e^{A^\top(\t-s)}q(\t)d\t,\eta(s)
\ran ds=-\i,$$
which contradicts \rf{4.16}. \endpf

\ms

Next, we would like to change the angle to look at the problem.

\bt{non-existence1} \sl Let $[A,B]$ be controllable and $q(\cd)\in L^1_{loc}(0,\i;\dbR^n)\setminus L^1(0,\i;\dbR^n)$ admit the polar decomposition
\bel{polar}q(s)=|q(s)|\th(s),\qq\forall s\in[t,\i),\ee
with $\th:[t,\i)\to S^n\equiv\{x\in\dbR^n\bigm||x|=1\}$. Suppose there exists an $\eta\in\dbR^n$ with
\bel{Aeta}A\eta\in\sR(B),\ee
and $\e>0$ such that
\bel{Ge<i}\lim_{T\to\i}\(\e\int_{G_\e\cap[0,T]}|q(s)|ds-|\eta|\int_{G_\e^c\cap[0,T]}
|q(s)|ds\)=\i,\ee
where
\bel{Ge}G_\e=\big\{s\in[0,\i)\bigm|\lan\th(s),\eta\ran\ges\e\big\}.\ee
Further, for some $\m>0$,
\bel{-m}\int_t^\i e^{-\m s}|q(s)|ds<\i.\ee
Then Problem (LQ) does not have an overtaking optimal control over $\sU_0[t,\i)$.

\et

\it Proof. \rm Suppose Problem (LQ) admits an overtaking optimal pair $(\bar X(\cd),\bar u(\cd))$. By the controllability of $[A,B]$, there exists a $\Th\in\dbR^{m\times n}$ such that $A_{\Th}\equiv A+B\Th$ is stable with
$$|e^{A_\Th s}|\les Me^{-\m s},\qq s\ges t,$$
for some $M\ges1$ and $\m>0$ (so that \rf{-m} holds). Define
$$\bar v(s)=\bar u(s)-\Th\bar X(s),\qq s\in[t,\i).$$
Then $(\bar X(\cd),\bar v(\cd))$ satisfies
$$\dot{\bar X}(s)=A_\Th\bar X(s)+B\bar v(s),\qq s\ges t,$$
and with
$$u(s)=\Th X(s)+v(s),\qq s\ges t,$$
we have
$$\ba{ll}
\ns\ds J_{_T}(t,x;\bar u(\cd))-J_{_T}(t,x;u(\cd))\\
\ns\ds=\int_t^T\(\lan Q\bar X(s),\bar X(s)\ran+2\lan S\bar X(s),\Th\bar X(s)+\bar v(s)\ran+|\Th\bar X(s)+\bar v(s)|^2+2\lan q(s),\bar X(s)\ran\)ds\\
\ns\ds\qq-\int_t^T\(\lan QX(s),X(s)\ran+2\lan SX(s),\Th X(s)+v(s)\ran+|\Th X(s)+v(s)|^2+2\lan q(s),X(s)\ran\)ds\\
\ns\ds\equiv\int_t^T\(g^\Th(s,\bar X(s),\bar v(s))-g^\Th(s,X(s),
v(s))\)ds\equiv J^\Th_{_T}(t,x;\bar v(\cd))-J^\Th_{_T}(t,x;v(\cd)),\ea$$
where
$$\ba{ll}
\ns\ds g^\Th(s,X,v)=\lan Q_\Th X,X\ran+2\lan S_\Th X,v\ran+|v|^2+2
\lan q(s),X\ran\\
\ns\ds\qq\qq\q\equiv\lan(Q+\Th^\top S+S^\top\Th+\Th^\top\Th)X,X\ran+2\lan(S+\Th)X,v\ran+|v|^2
+2\lan q(s),X\ran.\ea$$
Note that
$$Q_\Th-S_\Th^\top S_\Th=Q-S^\top S.$$
Hence, $(\bar X(\cd),\bar v(\cd))$ is an overtaking optimal control of Problem (LQ) with the system $[A_\Th,B]$ and the running cost rate function $g^\Th(s,X,v)$. Next, for some $\h v_0\in\dbR^m$,
$$A\eta=B\h v_0.$$
Let $v_0=\h v_0+\Th\eta$ which leads to $\eta=A_\Th^{-1}Bv_0$. Thus,
$$\lan\th(s),A_\Th^{-1}Bv_0\ran=\lan\th(s),\eta\ran\ges\e,\qq s\in G_\e.$$
Now, we take
$$u(\t)=\bar u(\t)-v_0{\bf1}_{[t,T]}(\t),\qq\t\in[t,\i).$$
Then $u(\cd)\in\sU_0[t,\i)$ and
$$\ba{ll}
\ns\ds\xi(s)=-\int_t^se^{A_\Th(s-\t)}Bv_0d\t=-e^{A_\Th s}\[\int_t^se^{-A_\Th\t}d\t\]Bv_0
=-e^{A_\Th s}\big[e^{-A_\Th t}-e^{-A_\Th s}\big]A_\Th^{-1}Bv_0\\
\ns\ds\qq=\big[I-e^{A_\Th(s-t)}\big]A_\Th^{-1}Bv_0,\qq s\in[t,T].\ea$$
Consequently,
$$\lan q(s),\xi(s)\ran=|q(s)|\lan\th(s),\big[I-e^{A_\Th(s-t)}\big]A_\Th^{-1}Bv_0\ran\ges|q(s)|\lan\th(s),\eta\ran-M|\eta|
e^{-\m(s-t)}|q(s)|.$$
Then, it follows that (making use of \rf{Ge<i}, and noting the integrals over $[0,t]$ are finite)
$$\int_t^T\lan q(s),\xi(s)\ran ds\ges\e\int_{G_\e\cap[t,T]}|q(s)|ds-|\eta|\int_{G_\e^c
\cap[t,T]}|q(s)|ds-M|\eta|e^{\m t}\int_t^\i
e^{-\m s}|q(s)|ds\to\i.$$
Hence, by taking $T>t$ large, we obtain that the right-hand side of the above is arbitrarily large, contradicting \rf{limJ-J<0}. Therefore, there is no overtaking optimal control over $\sU_0[t,\i)$. \endpf

\ms

Note that for any $\eta\in\dbR^n$, either $\lan\th(s),\eta\ran\ges0$, or $\lan\th(s),\eta\ran\les0$. Condition \rf{Ge<i} roughly means that either $|q(s)|$ grows uniformly fast and the directions $\th(s)$ are not uniformly distributed, as $s\to\i$, or although $\th(s)$ are almost uniformly distributed, but for the directions that $\lan\th(s),\eta\ran>0$, $|q(s)|$ grows faster.

\bc{} \sl Let $[A,B]$ be controllable and $q(\cd)\in L^1_{loc}(0,\i;\dbR^n)$ such that \rf{polar} holds with
\bel{lim th}\lim_{s\to\i}\th(s)=\th_0,\qq A\th_0\in\sR(B).\ee
Suppose
\bel{int|q|}\int_0^\i|q(s)|ds=\i,\qq\int_0^\i e^{-\m s}|q(s)|ds<\i,\ee
for some $\m>0$. Then Problem (LQ) does not have an overtaking optimal control over $\sU_0[t,\i)$.

\ec

\it Proof. \rm It suffices to take $\eta=\th_0$ in the above theorem. \endpf

\ms

The above results show that $\sR(B)$ plays an important role. The following gives some further result relevant to this.

\bt{} \sl Let {\rm(H)} hold. Let $q(\cd)\in L^1_{loc}(0,\i;\dbR^n)\setminus L^1(0,\i;\dbR^n)$ admit the polar decomposition \rf{polar}. Suppose $\th(\cd)$ is differentiable such that for some $r_0>0$,
\bel{th'-Ath}\dot\th(s)-A\th(s)\in B\big(\bar\cO(0,r_0)\big),\qq s\in[0,\i),\ee
where $\bar\cO(0,r_0))=\{v\in\dbR^m\bigm||v|\les r_0\}$. Further, there exists a $\d>0$ such that
\bel{lim}\lim_{T\to\i}{\int_T^{T+\d}|q(s)|ds\over\int_t^T|q(s)|ds}=0.\ee
Then Problem (LQ) does not have an overtaking optimal control over $\sU_0[t,\i)$.

\et

\it Proof. \rm Let $v:[t,\i)\to\dbR^m$ such that
$$\dot\th(s)-A\th(s)=Bv(s),\qq s\in[t,\i).$$
Define
$$u(\cd)=\bar u(\cd)+v(\cd){\bf1}_{[t,T]}(\cd)+\h v(\cd){\bf1}_{[T,T+\d)}(\cd),$$
where
$$\h v(\t)=-B^\top e^{A^\top(T+\d-\t)}W(\d)^{-1}[\th(T)-e^{A(T-t)}\th(t)],\qq W(\d)=\int_0^\d e^{A(\d-\t)}BB^\top e^{A^\top(\d-\t)}d\t.$$
By the controllability of $[A,B]$, we have the invertibility of $W(\d)$. Then
$$\ba{ll}
\ns\ds\xi(s)=\int_t^se^{A(s-\t)}Bv(\t)d\t=\int_t^se^{A(s-\t)}\big[\dot\th(\t)-A\th(\t)
\big]d\t\\
\ns\ds\qq=\int_t^se^{A(s-\t)}d\th(\t)-A\int_t^se^{A(s-\t)}\th(\t)d\t=\th(s)-e^{A(s-t)}\th(t),\q s\in[t,T],\ea$$
and
$$\ba{ll}
\ns\ds\xi(s)=e^{A(s-T)}\xi(T)-\int_T^se^{A(s-\t)}BB^\top e^{A^\top(T+\d-\t)}W(\d)^{-1}[\th(T)-e^{A(T-t)}\th(t)]d\t,\q s\in[T,T+\d),\\
\ns\ds\xi(s)=0,\qq\forall s\in[T+\d,\i).\ea$$
Hence, for any $\h T>T+\delta$, one has
$$\ba{ll}
\ns\ds\int_t^{\h T}\lan q(s),\xi(s)\ran ds=\int_t^T|q(s)|\lan\th(s),\xi(s)\ran ds+\int_T^{T+\d}|q(s)|\lan\th(s),\xi(s)\ran ds\\
\ns\ds=\int_t^T|q(s)|\big(1-\lan\th(s),e^{A(s-t)}\th(t)\ran\big) ds+\int_T^{T+\d}|q(s)|\lan\th(s),e^{A(s-T)}\xi(T)\ran ds
\\
\ns\ds\qq -\int_T^{T+\d}|q(s)|\lan\th(s),\int_T^se^{A(s-\t)}BB^\top e^{A^\top(T+\d-\t)}
W(\d)^{-1}\xi(T)d\t\ran ds\\
\ns\ds\ges\int_t^T|q(s)|(1-Me^{-\m(s-t)})ds-M\int_T^{T+\d}|q(s)||\th(T)-e^{A(T-t)}\th(t)|ds\\
\ns\ds\qq\q-\int_T^{T+\d}\(|q(s)|
\int_T^{T+\d}M^2\|B\|^2e^{-\m(s-\t)}e^{-\m(T+\d-\t)}d\t\)
\|W(\d)^{-1}\|\,|\th(T)-e^{A(T-t)}\th(t)|ds. \ea$$
Let $\k={1\over\m}\ln(4M)$. Then ${1\over4}-Me^{-\m\k}=0$, and $s\ges t+\k$ if and only if ${1\over4}\ges Me^{-\m(s-t)}$. Hence,
$$\ba{ll}
\ns\ds\int_t^T|q(s)|(1\1n-\1n Me^{-\m(s-t)})ds\1n=\1n{3\over4}
\int_t^T\2n|q(s)|ds\1n
+\2n\int_t^{t+\k}\2n|q(s)|\({1\over4}\1n-\1n Me^{-\m(s-t)}\)ds\1n+
\2n\int_{t+\k}^T\2n
|q(s)|\({1\over4}\1n-\1n Me^{-\m(s-t)}\)ds\\
\ns\ds\qq\qq\qq\qq\qq\qq\ges{3\over4}
\int_t^T\2n|q(s)|ds-\(M-{1\over4}\)\int_t^{t+\k}|q(s)|ds\ges
{1\over2}\int_t^T|q(s)|ds,\ea$$
provided $T$ is large. On the other hand,
$$\ba{ll}
\ns\ds \int_T^{T+\d}\(|q(s)|
\int_T^{T+\d}M^2\|B\|^2e^{-\m(s-\t)}e^{-\m(T+\d-\t)}d\t\)\|W(\d)^{-1}\|
|\th(T)-e^{A(T-t)}\th(t)|ds\\
\ns\ds\qq +M\int_T^{T+\d}|q(s)||\th(T)-e^{A(T-t)}\th(t)|ds\\
\ns\ds\les \[M^2\|B\|^2\|W(\d)^{-1}\|(1+M)\d+M(1+M)\]\int_T^{T+\d}
|q(s)|ds\equiv K(\d)\int_T^{T+\d}|q(s)|ds.\ea$$
Consequently,
$$\ba{ll}
\ns\ds\int_t^{\h T}\lan q(s),\xi(s)\ran ds\ges{1\over2}\int_t^T|q(s)|ds
-K(\d)\int_T^{T+\d}|q(s)|ds\\
\ns\ds\qq\qq\qq\q\ ={1\over2}\(\int_t^T|q(s)|ds\)\(1-
K(\d){\int_T^{T+\d}|q(s)|ds\over\int_t^T|q(s)|ds}\)\to\i,\q \hb{as}\ T\rightarrow\i.\ea$$
Hence, there exists a large enough $T_0>0$ such that
\bel{J-J>0*}\ba{ll}
\ns\ds J_{_{\h T}}(t,x;\bar u(\cd))-J_{_{\h T}}\big(t,x;\bar u(\cd)+\n[u(\cd)-\bar u(\cd)]\big)\\
\ns\ds\ges2\int_t^{\h T}\lan q(s),\xi(s)\ran ds-C\(|x|^2+\int_t^\i|\bar u(s)|^2ds\)^{1\over2}\(\int_t^\i|u(s)-\bar u(s)|^2ds\)^{1\over2}\\

\ns\ds\ges\int_t^{\h T}\lan q(s),\xi(s)\ran ds-C>0,\qq\forall\h T>T_0.\ea\ee
This is a contradiction. \endpf

\ms

Observe that condition \rf{lim} holds for functions of polynomial growth. For example, if
$$\lim_{s\to\i}{|q(s)|\over s^k}=q_0,$$
for some $q_0>0$, then for $T>t$ large enough,
$$\int_T^{T+\d}|q(s)|ds\les 2q_0\int_T^{T+\d}s^kds={2q_0\over k+1}\big[(T+\d)^{k+1}
-T^{k+1}\big]={2q_0\over k+1}\big[(k+1)\d T^k+\cds\big]\les CT^k,$$
and
$$\int_t^T|q(s)|ds\ges{q_0\over2}\int_{T\over2}^Ts^kds={q_0\over2(k+1)}\big[1-2^{-(k+1)}\big]
T^{k+1}.$$
Hence, \rf{lim} holds. However, if $|q(s)|$ grows exponentially, \rf{lim} may fail. For example, if $|q(s)|=e^s$, then
$$\lim_{T\to\i}{\int_T^{T+\d}|q(s)|ds\over\int_t^T|q(s)|ds}=\lim_{T\to\i}{e^{T+\d}-e^T\over e^T-e^t}\to e^\d-1>0.$$
Let us now explore some other cases. Suppose our control is more powerful in the system, by which we mean that $B$ is invertible. For such a case, without loss of generality, we assume that $m=n$ and $B=I$.

\ms

Note that the above proof also works for the case that $C([t,T];\dbR^n)$ is replaced by $C_u([t,\i);\dbR^n)$ of bounded uniformly continuous functions. We now present the following result.

\bt{non-existsnce2} \sl Let {\rm(H)} hold with $B=I$. Let $q(\cd)\in L^1_{loc}(0,\i;\dbR^n)$ admit polar decomposition \rf{polar}. Suppose $\th(\cd)$ is uniformly continuous and
\bel{|q|=i}\liminf_{T\to\infty}{1\over T}\int_0^T|q(s)|ds=\i.\ee
Then Problem (LQ) does not have an overtaking optimal control over $\sU_0[t,\i)$.

\ms

\et

\it Proof. \rm  Suppose Problem (LQ) admits an overtaking optimal control $\bar u(\cd)$ over $\sU_0[t,\i)$ with $\bar X(\cd)$ being the corresponding state trajectory. Since $\th(\cd)$ is uniformly continuous, for any $\e>0$, there exists $\d=\d(\e)>0$ such that
\bel{th-th}|\th(s)-\th(s')|<\e,\qq\forall s,s'\in[t,\i),~|s-s'|\les\d.\ee
Denote $t_i=t+i\d$, $i\ges0$. Then we let
$$W(\d)=\int_{t_i}^{t_{i+1}}e^{A(t_{i+1}-\t)}e^{A^\top(t_{i+1}-\t)}d\t=\int_0^\d e^{A\t}e^{A^\top\t}d\t,$$
which is invertible. Moreover, we note that
$$\|e^{A^\top\t}-I\|\les\sum_{k=1}^\i{\|A\|^k\t^k\over k!}\les\|A\|\t e^{\|A\|\t}\les\|A\|e^{\|A\|\d}\d,\qq\t\in[0,\d].$$
Hence, for any $x\in\dbR^n$,
$$\ba{ll}
\ns\ds\lan W(\d)x,x\ran=\int_0^\d|e^{A^\top\t}x|^2d\t=\int_0^\d\(|x|^2+2\lan(e^{A^\top\t}-I)x,x\ran
+|(e^{A^\top\t}-I)x|^2\)d\t\\
\ns\ds\ges\int_0^\d\(1-2\|e^{A^\top\t}-I\|-\|e^{A^\top\t}-I\|^2\)|x|^2d\t\ges\(1-2\|A\|e^{\|A\|\d}\d-\|A\|^2e^{2\|A\|\d}\d^2\)\d|x|^2\ges{\d\over2}
|x|^2,\ea$$
provided $\d>0$ is small. Then, by the symmetry of $W(\d)$, one has
\bel{W^{-1}}\|W(\d)^{-1}\|\les{2\over\d}.\ee
Now, for any $T>t$, let $\d>0$ be small and for some integer $N$, $T-t=N\d$. We define
$$\ba{ll}
\ns\ds u(\t)=\bar u(\t)+e^{A^\top(t_1-\t)}W(\d)^{-1}\th(t_1){\bf1}_{[t_0,t_1)}(\t)+\sum_{i=2}^{N-1}e^{A^\top(t_i-\t)}W(\d)^{-1}
\big[\th(t_i)-e^{A\d}\th(t_{i-1})\big]{\bf1}_{[t_{i-1},t_i)}(\t)\\
\ns\ds\qq\q-e^{A^\top(t_N-\t)}W(\d)^{-1}e^{A\d}\th(_{N-1}){\bf1}_{[t_{N-1},t_N)}(\t),\qq\t\in[t_0,\i).\ea$$
We let $v(\cd)=u(\cd)-\bar u(\cd)$ and
$$\xi(s)=\int_t^se^{A(s-\t)}v(\t)d\t,\qq s\in[t,\i).$$
We claim that
\bel{xi=th}\xi(t_k)=\th(t_k),\q\forall k=1,2,\cds,N-1;\qq\xi(s)=0,\q s\in[T,\i).\ee
In fact, the claim for $s\in[T,\i)$ is clear. For $s\in[t_0,t_1]$, we have
$$\ba{ll}
\ns\ds\xi(s)\equiv\int_{t_0}^se^{A(s-\t)}v(\t)d\t=\int_{t_0}^se^{A(s-\t)}e^{A^\top
(t_1-\t)}W(\d)^{-1}\th(t_1)d\t,\ea$$
which leads to $\xi(t_1)=\th(t_1)$. Next, for $s\in[t_1,t_2]$, we have
$$\xi(s)=e^{As}\xi(t_1)+\int_{t_1}^se^{A(s-\t)}v(\t)d\t=e^{As}\th(t_1)+\int_{t_1}^se^{A(s-\t)} e^{A^\top(t_2-\t)}W(\d)^{-1}[\th(t_2)-e^{A\d}\th(t_1)]d\t,$$
which leads to $\xi(t_2)=\th(t_2)$. By induction, suppose $\xi(t_k)=\th(t_k)$. Then for $s\in[t_k,t_{k+1}]$, one has
$$\ba{ll}
\ns\ds\xi(s)=e^{A(s-t_k)}\xi(t_k)+\int_{t_k}^se^{A(s-\t)}Bv(\t)d\t\\
\ns\ds\qq=e^{A(s-t_k)}\th(t_k)+\int_{t_k}^se^{A(s-\t)}e^{A^\top
(t_{k+1}-\t)}W(\d)^{-1}[\th(t_{k+1})-e^{A\d}\th(t_k)]d\t,\ea$$
which yields $\xi(t_{k+1})=\th(t_{k+1})$. This completes the induction. Consequently, for any $s\in[t,T]$, suppose $s\in[t_k,t_{k+1})$, one has
\bel{4.26}\ba{ll}
\ns\ds|\th(s)-\xi(s)|=\big|\th(s)-e^{A(s-t_k)}\th(t_k)-\int_{t_k}^s
e^{A(s-\t)}e^{A^\top
(t_{k+1}-\t)}W(\d)^{-1}[\th(t_{k+1})-e^{A\d}\th(t_k)]d\t\big|\\
\ns\ds\les|\th(s)-\th(t_k)|+\|I-e^{A(s-t_k)}\|+{2M^2\over\d}\(\int_{t_k}^se^{-\m(s-\t)}
e^{-\m(t_{k+1}-\t)}d\t\)|\th(t_{k+1})-e^{A\d}\th(t_k)|\\
\ns\ds\les|\th(s)\1n-\1n\th(t_k)|\1n+\1n\|A\|e^{\|A\|\d}\d\1n+\1n2M^2\(\|A\|e^{\|A\|\d}\d\1n
+\1n|\th(t_{k+1})\1n-
\1n\th(t_k)|\)\1n<\1n(1\1n+\1n2M^2)\(\e\1n+\1n\|A\|e^{\|A\|\d}\d\).\ea\ee
Consequently, for any $s\in[t,T]$, suppose $s\in[t_k,t_{k+1})$, one has
\bel{4.27}\ba{ll}
\ns\ds\lan q(s),\xi(s)\ran=|q(s)|\lan\th(s),\th(s)+\xi(s)-\th(s)\ran\ges|q(s)|\[1-(1+2M^2)\(\e+\|A|e^{\|A\|\d}\d\)\]\ges{1\over2}|q(s)|,\ea\ee
provided $\e,\d>0$ small enough. Note that
\bel{|v|}\ba{ll}
\ns\ds|v(\t)|^2\1n\les\1n\|e^{A^\top\1n(t_1\1n-\1n\t)}W(\d)^{-1}\|^2{\bf1}_{[t_0,t_1)}(\t)\1n
+\2n\sum_{i=2}^{N-1}\|e^{A^\top(t_i-\t)}W(\d)^{-1}\|^2|\th(t_i)-e^{A\d}
\th(t_{i-1})|^2{\bf1}_{[t_{i-1},t_i)}(\t)\\
\ns\ds\qq\qq+\|e^{A^\top(t_N-\t)}W(\d)^{-1}e^{A\d}\|^2{\bf1}_{[t_{N-1},t_N)}(\t)\\
\ns\ds\qq\q\les\1n {4M^2\over\d^2}\1n\({\bf1}_{[t_0,t_1)}(\t)\1n+\2n\sum_{i=2}^{N-1}\1n|\th(t_i)\1n-\1n e^{A\d}\th(t_{i-1})|^2
{\bf1}_{[t_{i-1},t_i)}(\t)\1n+\1n\|e^{A\d}\|^2{\bf1}_{[t_{N-1},t_N)}(\t)\),\q\t\in[t,\i).\ea\ee
Thus, $v(\cd)\in\sU^\i[t,\i)\cap\sU_0[t,\i)$ with (note $N\d=T-t$)
$$\int_t^T\3n|v(\t)|^2d\t\1n\les\1n {4M^2\over\d^2}\(\d\1n+\1n\d\sum_{i=2}^{N-1}|\th(t_i)\1n-\1n e^{A\d}\th(t_{i-1})|^2\1n+\1n\d\|e^{A\d}\|^2\)    \1n\les\1n{4M^2\over\d}\[1\1n+\1n\|e^{A\d}\|^2\1n+\1n{\big(\e+\|A\|e^{\|A\|\d}\d\big)^2\over
\d}(T\1n-\1n t)\].$$
Then take $\n=-1$ in \rf{J-J*}, we have
$$\ba{ll}
\ns\ds J_T(t,x;\bar u(\cd))-J_T\big(t,x;\bar u(\cd)-[u(\cd)-\bar u(\cd)]\big)\\
\ns\ds\ges2\int_t^T\lan q(s),\xi(s)\ran ds-C\(|x|^2+\int_t^T|\bar u(s)|^2ds\)^{1\over2}\(1+\int_t^T|u(s)-\bar u(s)|^2ds\)\\
\ns\ds\ges\1n\int_t^T\2n|q(s)|ds\1n-\1n C\[|x|\1n+\1n\|\bar u (\cd)\|_2\]\(1\1n+\1n\int_t^T\2n|v(s)|^2ds\)\1n\ges\1n\int_t^T\2n|q(s)|ds\1n-\1n K(\e,\d)\(|x|\1n+\1n\|\bar u(\cd)\|_2\)\big[1\1n+\1n(T\1n-\1n t)\big].\ea$$
Here, $K(\e,\d)$ is a constant independent of $T>t$. Thus, by our condition, for large enough $T>t$, the above will be larger than 0, leading to a contradiction. \endpf

\ms

Note that in the case $B=I$, $[A,B]$ is always controllable. Now, instead, if $m<n$, and $[A,B]$ is controllable, we do not have estimate \rf{W^{-1}}. Consequently, \rf{4.26} might not be true, and \rf{4.27} cannot be obtained this way. Therefore, we do not have the conclusion of the theorem.

\ms

\section{Existence of Overtaking Optimal Controls.}

From the previous section, we have seen that when $q(\cd)\in L^1_{loc}(0,\i;\dbR^n)\setminus L^1(0,\i;\dbR^n)$, the corresponding Problem (LQ) associated with \rf{state3}--\rf{I-SS>0,q} might not have overtaking optimal control over subspace $\sU[t,\i)$ of $\sU_0[t,\i)$ in general. In this section, we present a situation that the control set $\sU[t,\i)$ is a convex and closed subset of $\sU_0[t,\i)$, and Problem (LQ) has an overtaking optimal control.

\ms

Let \rf{e^A} hold and $q(\cd)$ grow at most exponentially with a rate $\a\in(0,\m)$. Then
\bel{q<i}\int_s^\i|e^{A^\top\t}q(\t)|d\t<\i.\ee
For such a case, we may let
\bel{hrho}\h\rho(s)=B^\top e^{-A^\top s}\int_s^\i e^{A^\top\t}q(\t)d\t,\qq s\in[t,\i),\qq\int_t^\i|\h\rho(s)|ds=\i.\ee

\ms

From Theorem \ref{Th-4.4}, we see that if $\bar u(\cd)$ is an interior point of $\sU[t,\i)$, most likely, it will not be overtaking optimal for our Problem (LQ). Now, for convex set $\sU[t,\i)\subseteq\sU_0[t,\i)$, we denote its boundary by $\pa\sU[t,\i)$. For any $u(\cd)\in\pa\sU[t,\i)$, an {\it extended inner normal} of $\sU[t,\i)$ at $u(\cd)$ is an element $\Bv(\cd)\in\sU_{loc}[t,\i)$ such that
\bel{Bv}\liminf_{T\to\i}\int_t^T\lan\Bv(s),v(s)-u(s)\ran ds\ges0,\qq\forall v(\cd)\in\sU[t,\i).\ee
We emphasize the above $\Bv(\cd)$. We let $\BN\big(u(\cd);\sU[t,\i)\big)\subseteq\sU_{loc}[t,\i)$ be the set of all extended inner normals of $\sU[t,\i)$ at $u(\cd)\in\pa\sU[t,\i)$. If $\Bv(\cd)\in\sU_0[t,\i)$ satisfies \rf{Bv}, then \rf{Bv} can be written as
\bel{Bv*}\int_t^\i\lan\Bv(s),v(s)-u(s)\ran ds\ges0,\qq\forall v(\cd)\in\sU[t,\i).\ee
Such a $\Bv(\cd)$ is called an {\it inner normal} of $\sU[t,\i)$ at $u(\cd)$. The set of all inner normals of $\sU[t,\i)$ at $u(\cd)\in\pa\sU[t,\i)$ is denoted by $\BN_0\big(u(\cd);\sU[t,\i)\big)$ which is a subset of $\sU_0[t,\i)$. Since in our LQ problem, $q(\cd)\in L^1_{loc}(0,\i;\dbR^n)\setminus L^1(0,\i;\dbR^n)$, it is easy for us to understand that (usual) inner normals are not enough and we need extended ones. On the other hand, let us recall that
$$\sU[t,T]=\big\{u(\cd)\big|_{[t,T]}\bigm|u(\cd)\in\sU[t,\i)\big\},$$
which it is convex and closed in $\sU_0[t,T]$. If $u(\cd)\in\pa\sU[t,\i)$, one must have $u(\cd)\big|_{[t,T]}\in\pa\sU[t,T]$. For any $T>t$, by taking $v(\cd)$ in such a way that it coincides with $u(\cd)$ on $[T,\i)$, we see that
\bel{Bv*}\int_t^T\lan\Bv(s),v(s)-u(s)\ran ds\ges0,\qq\forall v(\cd)\in\sU[t,T].\ee
Any $\Bv(\cd)\in\sU_0[t,T]$ satisfying \rf{Bv*} is called an {\it inner normal} of $\sU[t,T]$ at $u(\cd)\in\pa\sU[t,T]$. We let $\BN\big(u(\cd);\sU[t,T]\big)$ be the set of all inner normals of $\sU[t,T]$ at $u(\cd)$. Clearly, for any $u(\cd)\in\pa\sU[t,\i)$,
$$\BN\big(u(\cd);\sU[t,\i)\big)\big|_{[t,T]}\equiv\Big\{\Bv(\cd)\big|_{[t,T]}\bigm|
\Bv(\cd)\in\BN\big(u(\cd);\sU[t,\i)\big)\Big\}=\BN(u(\cd);\sU[t,T]\big).$$
But
$$\BN\big(u(\cd);\sU[t,\i)\big)=\bigcup_{T>t}\BN\big(u(\cd);\sU[t,T]\big)\ne
\BN_0(u(\cd);\sU[t,\i)\big).$$
Next, let us recall $F_0(s)x$ and $F_1[u(\cd)](s)$ from \rf{F}.
Note
$$\ba{ll}
\ns\ds F_1[\bar u(\cd)](s)\1n=\1n S\2n\int_t^s\2n e^{A(s-\t)}B\bar u(\t)d\t
\1n+\2n\int_s^\i\2n B^\top e^{A^\top(\t-s)}Q\2n\int_t^\t\2n e^{A(\t-r)}B\bar u(r)drd\t\1n+\2n\int_s^\i\2n B^\top e^{A^\top(\t-s)}S^\top\bar u(\t)d\t\\
\ns\ds=\int_t^sSe^{A(s-\t)}B\bar u(\t)d\t
+\int_t^\i\int_{s\vee\t}^\i B^\top e^{A^\top(r-s)}Qe^{A(r-\t)}B\bar u(\t)drd\t+\int_s^\i B^\top e^{A^\top(\t-s)}S^\top\bar u(\t)d\t\\
\ns\ds=\2n\int_t^\i\2n\[{\bf1}_{[t,s]}
(\t)Se^{A(s-\t)}\1n B\1n
+\2n{\bf1}_{[s,\i)}(\t)B^\top\1n e^{A^\top\1n(\t-s)}S^\top
\3n+\2n\int_{s\vee\t}^\i\2n B^\top\1n e^{A^\top\1n(r-s)}Qe^{A(r-\t)}Bdr\]\bar u(\t)d\t\1n\equiv\2n\int_t^\i\3n\F(s,\t)\bar u(\t)d\t.\ea$$
From \rf{F*}, one has
$$\int_t^\i\2n\Big|\1n\int_t^\i\2n\F(s,\t)\bar u(\t)d\t\Big|^2\1n ds\1n=\1n\int_t^T\1n\big|F_1[\bar u(\cd)](s)\big|^2ds\1n\les\1n{3\|B\|^2M^2\over\m^2}\({\|B\|^2M^2\|Q\|^2\over\m^2}\1n
+\1n2\|S\|^2\)\1n
\int_t^\i\2n|\bar u(s)|^2ds.$$
Thus, by choosing $\m>0$ in \rf{e^A} large enough, we may assume that
\bel{k}\k={3\|B\|^2M^2\over\m^2}\({\|B\|^2M^2\|Q\|^2\over\m^2}+2\|S\|^2\)\in(0,1).\ee
Hence,
$$\int_t^\i\Big|\int_t^\i\F(s,\t)\bar u(\t)d\t\Big|^2ds\les\k
\int_t^\i|\bar u(s)|^2ds,\qq\forall\bar u(\cd)\in\sU_0[t,\i),$$
for some $\k\in(0,1)$. This leads to that for any $\f(\cd)\in\sU_0[t,\i)$, the following Fredholm integral equation of the second kind
$$\f(s)+u(s)+\int_t^\i\F(s,\t)u(\t)d\t=0$$
admits a unique solution $u(\cd)\in\sU_0[t,\i)$. We now state the main theorem of this section.

\bt{} \sl Let {\rm(H)}, \rf{q<i} and \rf{k} hold. Let $\h\rho(\cd)$ be defined by \rf{hrho}. Let $\bar u(\cd)\in\pa\sU[t,\i)$ be the solution to the following Fredholm integral equation of the second kind:
\bel{int}\h\rho_0(s)+F_0(s)x+\bar u(s)+\int_t^\i\F(s,\t)\bar u(\t)d\t=0,\qq s\in[t,\i),\ee
for some $\h\rho_0(\cd)\in\sU_0[t,\i)$ such that
\bel{hrho_1}\h\rho_1(\cd)\equiv\h\rho(\cd)-\h\rho_0(\cd)\in\BN(\bar u(\cd);\sU[t,\i)\big).\ee
Then $\bar u(\cd)$ is an overtaking optimal control of Problem (LQ) at $(t,x)\in[0,\i)\times\dbR^n$ over $\sU[t,\i)$.

\et

\it Proof. \rm By Theorem \ref{Th-4.2}, we know that $\bar u(\cd)\in\sU[t,\i)$ is overtaking optimal if and only if
\bel{5.11}0\les\liminf_{T\to\i}\int_t^T\lan\h\rho(s)+F_0(s)x+\bar u(s)+\int_t^\i\F(s,\t)\bar u(\t)d\t,u(s)-\bar u(s)\ran ds,\q\forall u(\cd)\in\sU[t,\i).\ee
From our condition, we have
$$\int_t^T\lan\h\rho(s)+F_0(s)x+\bar u(s)+\int_t^\i\F(s,\t)\bar u(\t)d\t,u(s)-\bar u(s)\ran ds=\int_t^T\lan\h\rho_1(s),u(s)-\bar u(s)\ran ds.$$
Hence, \rf{5.11} holds leading to that $\bar u(\cd)$ is overtaking optimal. \endpf

\ms

It is not hard to see that the conditions assumed in the above theorem can be replaced by the following: There exists a $T_0>t$ such that
\bel{}\h\rho(s){\bf1}_{[t,T_0]}(s)+F_0(s)x+\bar u(s)+\int_t^\i\F(s,\t)\bar u(\t)d\t=0,\qq s\in[t,\i),\ee
admits a unique solution $\bar u(\cd)\in\pa\sU[t,\i)$, and
\bel{}\h\rho(s)\in\BN\big(\bar u(s);\sU[t,\i)\big),\qq s\in[T_0,\i).\ee
Then $\bar u(\cd)$ is an overtaking optimal control for Problem (LQ) at
$(t,x)\in[0,\i)\times\dbR^n$ over $\sU[t,\i)$.

Thus, we would like to have
$$0\les\liminf_{T\to\i}\int_t^T\lan\h\rho(s)+F_0(s)x+\bar u(s)
+\int_t^\i\F(s,\t)\bar u(\t)d\t,u(s)-\bar u(s)\ran ds,\qq\forall
u(\cd)\in\sU[t,\i).$$

Let us now present an illustrative example for the above theorem.

\bex{} \rm Consider controlled system in $\dbR^2$:
$$\left\{\2n\ba{ll}
\ds\dot X(s)=-\m X(s)+u(s),\\ %
\ns\ds X(t)=x,\ea\right.$$
with $\m>1$. Then
$$X(s)=e^{-\m(s-t)}x+\int_t^se^{-\m(s-\t)}u(\t)d\t.$$
Let $Q=I$, $S=0$, $U=\bigm\{u=(u_1,u_2)\bigm|u_1\ges0\bigm\}$, and $q(s)=\begin{pmatrix}e^s\\ e^{-s}\end{pmatrix}$. Then
$$\ba{ll}
\ns\ds\h\rho(s)=e^{\m s}\int_s^\i e^{-\m\t}q(\t)d\t=e^{\m s}\int_s^\i
\begin{pmatrix}e^{-(\m-1)\t}\\ e^{-(\m+1)\t}\end{pmatrix}d\t=\begin{pmatrix}{e^s\over\m-1}\\ {e^{-s}\over\m+1}\end{pmatrix},\\
\ns\ds F_0(s)x=\int_s^\i e^{-\m(\t-s)}e^{-\m(\t-t)}xd\t={e^{-\m(s-t)}\over2\m}x,\\
\ns\ds F_1[\bar u(\cd)](s)=\int_s^\i e^{-\m(\t-s)}\bar X_0(\t)d\t=\int_s^\i e^{-\m(\t-s)}\int_t^\t e^{-\m(\t-r)}\bar u(r)drd\t=\int_t^\i{e^{-\m|\t-s|}\over2\m}\bar u(\t)d\t.\ea$$
Let $\bar u(s)=(0,\bar u_2(\cd))\equiv\bar u_2(\cd)\Be_2$, with $\bar u_2(\cd)$ solves the following Fredholm integral equation:
$${e^{-s}\over\m+1}+{e^{-\m(s-t)}\over2\m}x_2+{e^{-s}\over\m+1}+\bar u_2(s)
+\int_t^\i{e^{-\m|s-\t|}\over2\m}\bar u(\t)d\t=0,\q s\in[t,\i).$$
Then, under such a $\bar u(\cd)$, we have
$$\ba{ll}
\ns\ds\liminf_{T\to\i}\int_t^T\lan\h\rho(s)+F_0(s)x+\bar u(s)+F_1[\bar u(\cd)](s),u(s)-\bar u(s)\ran ds\\
\ns\ds=\liminf_{T\to\i}\int_t^T\lan\begin{pmatrix}{e^s\over\m-1}+{e^{-\m(s-t)}\over2\m}x_1
\\ 0\end{pmatrix},u(s)-\bar u(s)\ran ds=\liminf_{T\to\i}\int_t^T\({e^s\over\m-1}+{e^{-\m(s-t)}\over2\m}x_1\)u_1(s)ds\ges0.\ea$$
This means that $\bar u(\cd)$ is an overtaking optimal control for Problem (LQ) at $(t,x)\in[0,\i)\times\dbR^2$ over $\sU[t,\i)\equiv\{u:[t,\i)\to U\bigm|
\int_t^\i|u(s)|^2ds<\i\}$.

\ex

\section{Concluding Remarks}

We have studied LQ problems over infinite time horizons for which no controllability/stabilizability condition is assumed for the homogeneous system $[A,B]$, no global integrability conditions are assumed for the nonhomogeneous term $b(\cd)$ in the state equation and the weighting functions $q(\cd),\rho(\cd)$ for the linear terms in the cost functional. For such kind of problems, all the existing approaches do not apply. We have established existence of overtaking optimal controls for the case that $(X_{_\Pi}(\cd),u(\cd))$ and $X_{\Pi^{\perp}}(\cd)$ are decoupled in the state equation and separated in the cost functional, and for a convex control set case. More interestingly, we also have proved the non-existence of overtaking optimal control for several cases. From these results, we see that overtaking optimality approach can be used to solve some problems that the classical tools are not applicable. However, the power of this approach is not unlimited.

{}


\vskip6mm
\begin{thebibliography}{99}

\bibitem{Anderson-Moore 1989} B.~D.~O.~Anderson and J.~B.~Moore, \sl Optimal control: Linear quadratic Methods, \rm Prentice-Hall International, Inc., 1989.

\bibitem{Arrow 1968} K.~J.~Arrow, \it Applications of control theory to economic growth, \sl Mathematics of the Decision Sciences, Part 2, \rm AMS, Providence, R.I., 1968, 85--119.

\bibitem{Arrow-Kurtz 1970} K.~J.~Arrow and M.~Kurtz, \it Optimal growth with irreversible investment in a Ramsey model, \sl Econometrica, \rm 38 (1970), 331--344.


\bibitem{Aseev-Kryazhimskii 2007} S.~M.~Aseev and A.~V.~Kryazhimskii, \it The Pontryagin maximum principle and problems of optimal economic growth, \sl Proc. Stekolv Inst. Math., \rm 257 (2007), 1--255.



\bibitem{Basco-Cannarsa-Frankowska 2018} V.~Basco, P.~Cannarsa, and H.~Frankowska, \it Necessary conditions for infinite horizon optimal control problems with state constraints, \sl Math. Control Relat. Fields, \rm 8 (2018), 535--555.

\bibitem{Belyakov 2019} A.~O.~Belyakov, \it On a sufficient condition for infinite horizon optimal control problems, \rm arXiv:1909.07379v1 [math.OC] 17 Sep 2019.

\bibitem{Brock-Haurie 1976} W.~A.~Brock and A.~Haurie, \it On existence of overtaking optimal trajectories over an infinite time horizon, \sl Math. Oper. Res., \rm 1 (1976), 337--346.

\bibitem{Buckdahn-Li-Quincampoix-Renault 2020} R.~Buckdahn, J.~Li, M.~Quincampoix, and J.~Renault, \it Representation formulas for limit values of long run stochastic optimal controls, \sl SIAM J. Control Optim., \rm 58 (2020), 1846--1873.


\bibitem{Carlson-Haurie-Leizarowitz 1991} D.~A.~Carlson, A.~B.~Haurie, and A.~Leizarowitz, \sl Infinite Horizon Optimal Control: Deterministic and Stochastic Systems, \rm Springer-Verlag, Berlin, 1991.





\bibitem{Halkin 1974} H.~Halkin, \it Necessary conditions for optimal control problems with infinite horizons, \sl Econometrica, \rm 42 (1974), 267--272.

\bibitem{Huang-Li-Yong 2015} J.~Huang, X.~Li, and J.~Yong, \it A linear-quadratic optimal control problem for mean-field stochastic differential equations in infinite horizon, \sl Math. Control Relat. Fields, \rm 5 (2015), 97--139.


\bibitem{Kalman 1960} R.~E.~Kalman, \it Contributions to the theory of optimal control, \sl Bol. Soc. Mat. Mexicana, \rm 5 (1960), 102--119.




\bibitem{Ramsey 1928} F.~P.~Ramsey, \it A mathematical theory of saving, \sl Economic J., \rm 38 (1928), 543--559.


\bibitem{Samuelson-Moussavian 1985} L.~Samuelson and M.~Moussavian, \it Sufficient conditions in optimal control theory with unbounded objectives, \sl Int. Econ. Review, \rm 26 (1985), 271--275.

\bibitem{Skritek-Veliov-2015} B.~Skritek and V.~M.~Veliov, \it On the infinite-horizon optimal control of age-structured systems, \sl J. Optim. Theory Appl., \rm 167    (2015), 243--271.

\bibitem{Sun-Yong 2018} J.~Sun and J.~Yong, \it Stochastic linear quadratic optimal control problems in infinite horizon, \sl Appl. Math. Optim., \rm 78 (2018), 145--183.

\bibitem{Sun-Yong 2020} J.~Sun and J.~Yong, \sl Stochastic Linear-Quadratic Optimal Control Theory: Open-Loop and Closed-Loop Solutions, \rm Springer, 2020.

\bibitem{Tan-Rugh 1998} H.~Tan, W.~J.~Rugh, \it On overtaking optimal tracking for linear systems, \sl Systems Control Lett., \rm 33 (1998), 63--72.

\bibitem{Weizsacker 1965} C.~C.~von Weizs\"acker, \it Existence of optimal progras of accumulation for an infinite time horizon, \sl Rev. Econom. Stud.,  \rm 32 (1965), 85--104.



\bibitem{Willems 1971} J.~C.~Willems, \it Least sqares stationary optimal control and the algebraic Riccati equation, \sl IEEE Trans. Automat. Control, \rm 16 (1971), 621--634.

\bibitem{Wonham 1985} W.~M.~Wonham, \sl Linear Multivariable Control: A Geometric Approach, 3rd Edition, \rm Springer Science+Business Media, LLC, 1985.




 \end{thebibliography}
\end{document}